\documentclass[11pt]{article}

\usepackage{amsmath,amssymb,amsthm}
\usepackage[margin=1in]{geometry}
\usepackage{booktabs}
\usepackage{listings}
\usepackage{enumerate}
\usepackage[hidelinks]{hyperref}

\newtheorem{theorem}{Theorem}[section]
\newtheorem{lemma}[theorem]{Lemma}
\newtheorem{proposition}[theorem]{Proposition}

\theoremstyle{definition}
\newtheorem{definition}[theorem]{Definition}
\newtheorem{remark}[theorem]{Remark}
\newtheorem{hypothesis}[theorem]{Hypothesis}

\theoremstyle{plain}
\newtheorem*{thmA}{Theorem A}
\newtheorem*{thmB}{Theorem B}
\newtheorem*{thmC}{Theorem C}
\newtheorem*{corD}{Corollary D}

\newcommand{\cl}{\mathrm{cl}}
\newcommand{\Irr}{\mathrm{Irr}}
\newcommand{\GF}{\mathrm{GF}}
\newcommand{\GL}{\mathrm{GL}}
\newcommand{\Gal}{\mathrm{Gal}}
\newcommand{\Sp}{\mathrm{Sp}}
\newcommand{\Ind}{\mathrm{Ind}}
\newcommand{\X}{\mathcal{X}}
\newcommand{\G}{\mathcal{G}}
\newcommand{\Ob}{\mathcal{O}}

\title{Large orbits of $p$-groups on conjugacy classes}
\author{Thomas Michael Keller\\[1mm]
\small Department of Mathematics, Texas State University\\
\small 601 University Drive, San Marcos, TX 78666, United States of America\\
\small \texttt{keller@txstate.edu}}
\date{}

\begin{document}

\maketitle

\begin{abstract}
Let $P$ be a $p$-group that acts faithfully and coprimely on a finite group $H$. A classical result of D.\ S.\ Passman, in the form proved by I.\ M.\ Isaacs, asserts that there exists an element $x\in H$ such that $|C_P(x)|\le (|P|/p)^{1/p}$. In this paper, we study the corresponding question for the action of $P$ on the conjugacy classes of $H$. We show that for every odd prime $p$ the analogue of Isaacs' bound fails for classes. On the other hand, we prove that if $H$ is solvable, then there always exists a conjugacy class $c$ of $H$ whose stabilizer in $P$ has order at most $|P|^{1/2}$, and that this bound is best possible for every prime $p$. As a consequence, there always exists an irreducible character of $H$ whose inertia group in $P$ has order at most $|P|^{1/2}$. The proof proceeds via a result on linear actions of $p$-nilpotent groups, and in the quasiprimitive case it relies on a classification of Holt and Yang and on computer calculations.
\end{abstract}


\section{Introduction}

One important aspect of the study of finite group actions is the existence of large orbits, or equivalently, of points with small stabilizers. A classical result of this type is due to D.\ S.\ Passman \cite{Pa}, who showed that if a $p$-group $P$ acts faithfully on a solvable $p'$-group $H$, then there exists $x\in H$ such that $|C_P(x)|\le|P|^{1/2}$ if $p=2$ and $|C_P(x)|\le |P|^{1/3}$ if $p>2$. At the end of his paper, Passman stated without proof that one always has $|C_P(x)|\le|P|^{1/p}$. This was proved by I.\ M.\ Isaacs \cite{Is99}, in fact in the following slightly stronger form, which moreover does not require $H$ to be solvable.

\begin{theorem}[Isaacs]\label{isaacs}
Let $P$ be a nontrivial $p$-group that acts faithfully on a group $H$, where $|H|$ is not divisible by $p$. Then there exists an element $x\in H$ such that $|C_P(x)|\le(|P|/p)^{1/p}$.
\end{theorem}

An essential ingredient in Isaacs' proof is a lemma of B.\ Hartley and A.\ Turull \cite{HT}, which allows one to replace $H$ by an abelian group without changing the permutation action of $P$ on the elements of $H$.

In this paper, we are interested in the action of $P$ on the set $\cl(H)$ of conjugacy classes of $H$ rather than on the elements of $H$. If $c\in\cl(H)$, we write $C_P(c)$ for the stabilizer of $c$ in $P$, that is, for the set of all $g\in P$ with $c^g=c$. Clearly $C_P(x)\le C_P(x^H)$ for every $x\in H$, so that a bound for the stabilizers of classes is a strengthening of a bound for the centralizers of elements. Note also that the Hartley--Turull lemma is of no help here, because replacing $H$ by an abelian group changes the classes completely.

In view of Theorem \ref{isaacs}, it is natural to ask whether there always exists a class $c$ of $H$ with $|C_P(c)|\le |P|^{1/p}$. Our first main result shows what the correct bound is, at least when $H$ is solvable.

\begin{thmA}
Let $P$ be a $p$-group that acts faithfully and coprimely on a finite solvable group $H$. Then there exists a conjugacy class $c$ of $H$ such that
$$|C_P(c)|\le|P|^{1/2}.$$
\end{thmA}

The following result shows that Theorem A is best possible for every prime $p$, and that for odd $p$ the bound of Theorem \ref{isaacs} does not carry over to classes.

\begin{thmB}
Let $p$ be a prime and let $q$ be a prime power such that $p$ divides $q-1$; if $p=2$, assume moreover that $q\equiv 1\pmod 4$. Let $F=\GF(q^p)$, let $U$ be the kernel of the norm map $F^\times\to\GF(q)^\times$, and let $K$ be the Hall $p'$-subgroup of the cyclic group $U$. Let $H=F\rtimes K$ be the group of affine maps $y\mapsto ky+w$ with $k\in K$ and $w\in F$, let $\omega\in\GF(q)^\times$ be of order $p$, and let $\varphi$ be the Frobenius map $y\mapsto y^q$ of $F$. Then $P=\langle\omega\rangle\times\langle\varphi\rangle\cong C_p\times C_p$ acts faithfully and coprimely on $H$ by conjugation, $H$ is solvable, and
$$|C_P(c)|\ge p=|P|^{1/2}\quad\text{for every } c\in\cl(H).$$
\end{thmB}

The smallest instance of Theorem B is $p=2$ and $q=5$, where $|H|=25\cdot3=75$. A particularly simple instance with $p$ odd and $q$ prime is $p=3$ and $q=7$. Here $|H|=7^3\cdot 19$ and $P\cong C_3\times C_3$, and each of the $37$ conjugacy classes of $H$ has a stabilizer of order at least $3>9^{1/3}$ in $P$, while, as predicted by Theorem \ref{isaacs}, $P$ has a regular orbit on the elements of $H$.

Theorem A will be obtained from a result on linear actions of $p$-nilpotent groups. In order to state it, we introduce the following class of pairs.

\begin{definition}\label{defX}
Let $p$ be a prime. We denote by $\X_p$ the class of all pairs $(G,V)$, where $G$ is a finite solvable $p$-nilpotent group with $O_p(G)>1$, and $V$ is a finite faithful irreducible $G$-module over a field $\GF(r)$ of prime order $r\neq p$.
\end{definition}

Note that if $G$ has an irreducible faithful module $V$ in characteristic $p$, then $O_p(G)$ acts trivially on $V$, so that $O_p(G)=1$. Thus the condition $r\neq p$ in Definition \ref{defX} is automatically satisfied. Also note that the orbit structure of $G$ on $V$ does not change if $V$ is regarded as a module over a smaller field, but such a module need not remain irreducible. For this reason we restrict ourselves to prime fields, which is all that is needed for Theorem A.

\begin{thmC}
Let $(G,V)\in\X_p$ and let $N$ be a normal subgroup of $G$. Then the following hold.
\begin{enumerate}[(a)]
\item There exists $v\in V$ such that $|C_N(v)|_p\le|N|_p^{1/2}$.
\item Suppose that $|N|_p=p^{2a}$ for some $a\ge 1$ and that $|C_N(v)|_p\ge p^a$ for all $v\in V$. Then the vectors $v\in V$ with $|C_N(v)|_p=p^a$ lie in at least two different $N$-orbits.
\end{enumerate}
\end{thmC}

Part (a) of Theorem C with $N=G$ is what is needed for Theorem A. Part (b) is needed to make the induction work, and we will see in Proposition \ref{wreath} that, at least for $N=G$, it is in fact a consequence of part (a) applied to suitable wreath products. Moreover, for prime $q$ the groups in Theorem B give rise to pairs $(G_0,V_0)\in\X_p$ such that $|C_{G_0}(v)|_p=|G_0|_p^{1/2}$ for all $0\neq v\in V_0$ (see Remark \ref{sharpC}). So the inequality in (a) cannot be replaced by a strict inequality.

The coprime actions of $P$ on the classes and on the irreducible characters of $H$ are known to be permutation isomorphic. Hence we also obtain the following.

\begin{corD}
Let $P$ be a $p$-group that acts faithfully and coprimely on a finite solvable group $H$. Then there exists $\chi\in\Irr(H)$ such that $|I_P(\chi)|\le|P|^{1/2}$, where $I_P(\chi)$ denotes the stabilizer of $\chi$ in $P$.
\end{corD}

Large orbits of coprime actions on characters were studied, for instance, by A.\ Moret\'o \cite{Mo} and in \cite{KY}. By Theorem B, the bound in Corollary D is best possible.

We now briefly describe the proof of Theorem C, which is by induction on $|G|$. If $V$ is not quasiprimitive, we show that either there is a normal subgroup of index $p$ in $G$ which stabilizes the members of a system of imprimitivity with $p$ components, or there is a system of imprimitivity whose stabilizers contain a Sylow $p$-subgroup acting faithfully on the corresponding component. In both cases we can glue together vectors obtained by induction, and here part (b) of Theorem C is essential. If $V$ is quasiprimitive and the invariant $e=|E/Z(E)|^{1/2}$ of the normal subgroup $E$ of extraspecial type in $G$ equals $1$, then $G$ is a subgroup of a semilinear group $\Gamma(V)$, and this case can be handled by an elementary counting argument. If $V$ is quasiprimitive and $e>1$, we use the classification of the quasiprimitive solvable linear groups without regular orbits due to D.\ F.\ Holt and Y.\ Yang \cite{HY}. Most of the exceptional cases in this classification are excluded by a short argument, and the remaining ones are handled by computer calculations with GAP \cite{GAP} and the IRREDSOL library \cite{IRREDSOL}. Since the classification in \cite{HY} excludes metacyclic groups, the metacyclic case is treated separately, again by a counting argument supplemented by computer calculations in a few small cases.

The paper is organized as follows. In Section 2 we collect some auxiliary results. In Section 3 we prove Theorem B. In Section 4 we reduce Theorem A to Theorem C. Sections 5 and 6 contain the proof of Theorem C in the imprimitive and quasiprimitive case, respectively, and in Section 7 we complete the proofs of our main results. We conclude with some remarks and open questions in Section 8.

\section{Auxiliary results}

Throughout, all groups and modules are finite. If $n$ is a positive integer and $p$ is a prime, we write $n_p$ for the $p$-part of $n$, and if $X$ is a group, we write $|X|_p$ for $|X|_p=(|X|)_p$. If a group $X$ acts on a set $\Omega$ and $\omega\in\Omega$, then $C_X(\omega)$ denotes the stabilizer of $\omega$ in $X$, and similarly for subsets of $\Omega$. If $V$ is a $G$-module and $N\le G$, then $V_N$ denotes the restriction of $V$ to $N$. A semisimple module is called homogeneous if all its irreducible constituents are isomorphic. Recall that by Clifford's theorem, if $V$ is irreducible and $N\trianglelefteq G$, then $V_N$ is semisimple, and its homogeneous components are permuted transitively by $G$. The module $V$ is called quasiprimitive if $V_N$ is homogeneous for all $N\trianglelefteq G$.

The following notation will be used throughout the paper.

\begin{definition}\label{tight}
Let $N$ be a group acting on a finite set $V$ (usually a vector space), and let $p$ be a prime. We put
$$\mu_p(N,V)=\min_{v\in V}|C_N(v)|_p,$$
and call $v\in V$ a minimizer if $|C_N(v)|_p=\mu_p(N,V)$. We say that the action of $N$ on $V$ is tight if $|N|_p=p^{2a}$ for some $a\ge 1$ and $\mu_p(N,V)=p^a$, and we say that it is rigid if it is tight and the minimizers form a single $N$-orbit.
\end{definition}

With this notation, Theorem C states that for $(G,V)\in\X_p$ and $N\trianglelefteq G$ we have $\mu_p(N,V)\le|N|_p^{1/2}$, and that the action of $N$ on $V$ is not rigid. Note that if $v\in V$ is such that $C_N(v)$ is a $p'$-group, then $\mu_p(N,V)=1$, and so both statements hold. In particular, Theorem C holds for $(G,V)$ if $G$ has a regular orbit on $V$.

The next lemma, which is elementary, connects the stabilizers of classes with centralizers in a semidirect product.

\begin{lemma}\label{classes}
Let $P$ be a $p$-group acting on a $p'$-group $H$, and let $G=HP$ be the corresponding semidirect product. Let $x\in H$ and let $c=x^H$ be the class of $x$ in $H$. Then
$$|C_P(c)|=|C_G(x)|_p.$$
\end{lemma}

\begin{proof}
Let $S$ be the stabilizer of $c$ in $G$, where $G$ acts on $\cl(H)$ by conjugation. Since $H$ fixes $c$ and $G=HP$, we have $S=H\,C_P(c)$, and hence $|S:H|=|C_P(c)|$. On the other hand, since $H$ is transitive on $c$, we have $S=C_G(x)H$. Thus
$$|C_P(c)|=|C_G(x)H:H|=|C_G(x):C_H(x)|.$$
As $C_H(x)$ is a $p'$-group and $C_G(x)/C_H(x)\cong C_G(x)H/H$ is isomorphic to a subgroup of the $p$-group $G/H$, we conclude that $|C_G(x):C_H(x)|=|C_G(x)|_p$, as desired.
\end{proof}

We will frequently use the following simple lemma, which is a variation of Lemma 2.1 of \cite{Is99}.

\begin{lemma}\label{chain}
Let a group $N$ act on a vector space $V=V_1\oplus\cdots\oplus V_m$, where the $V_j$ are $N$-invariant. Put $M_0=N$ and $M_j=C_N(V_1\oplus\cdots\oplus V_j)$ for $j=1,\dots,m$. Let $v_j\in V_j$ for all $j$ and $v=v_1+\cdots+v_m$. Then
$$|C_N(v)|_p\le|M_m|_p\prod_{j=1}^m|C_{M_{j-1}/M_j}(v_j)|_p.$$
\end{lemma}

\begin{proof}
For $j=0,\dots,m$ put $X_j=C_{M_j}(v_{j+1}+\cdots+v_m)$, so that $X_0=C_N(v)$ and $X_m=M_m$. Since the $V_j$ are $N$-invariant, an element of $N$ fixes $v$ if and only if it fixes every $v_j$. It follows that $X_{j-1}\cap M_j=X_j$ and that $X_{j-1}M_j/M_j\le C_{M_{j-1}/M_j}(v_j)$ for $j=1,\dots,m$. Hence
$$|X_{j-1}|_p=|X_{j-1}M_j/M_j|_p\,|X_j|_p\le|C_{M_{j-1}/M_j}(v_j)|_p\,|X_j|_p,$$
and the assertion follows by induction.
\end{proof}

\begin{lemma}\label{permmod}
Let $Z$ be a $p$-group, $Y<Z$, and let $W$ be a nonzero module for $Y$ over a field of characteristic different from $p$ on which $Y$ acts trivially. Then $\Ind_Y^Z(W)$ is not homogeneous.
\end{lemma}

\begin{proof}
The module $\Ind_Y^Z(W)$ is a direct sum of copies of the permutation module $M$ of $Z$ on the cosets of $Y$. By Maschke's theorem, $M$ is the direct sum of the trivial module and the augmentation submodule, which is nonzero since $|Z:Y|>1$. Since $Z$ acts nontrivially on the cosets of $Y$, it acts nontrivially on the augmentation submodule. Thus $\Ind_Y^Z(W)$ has both trivial and nontrivial irreducible constituents.
\end{proof}

If $K$ is a finite field and $\sigma$ is a field automorphism of $K$, then a map $y$ of a $K$-vector space $W$ is called $\sigma$-semilinear if it is additive and $y(\lambda w)=\sigma(\lambda)y(w)$ for all $\lambda\in K$ and $w\in W$.

\begin{lemma}\label{semilinear}
Let $K=\GF(Q)$, let $\sigma$ be an automorphism of $K$ of order $s$ with fixed field $K_0$, and let $y$ be a bijective $\sigma$-semilinear map of a $K$-vector space $W$ of dimension $m$. Then the set $C_W(y)$ of fixed points of $y$ is a $K_0$-subspace of $W$, and $|C_W(y)|\le Q^{m/s}$.
\end{lemma}

\begin{proof}
Clearly $C_W(y)$ is a $K_0$-subspace. We claim that elements of $C_W(y)$ which are linearly independent over $K_0$ are linearly independent over $K$. Otherwise, choose a nontrivial relation $\sum_{i=1}^t\lambda_iw_i=0$ with $w_i\in C_W(y)$ linearly independent over $K_0$, $\lambda_i\in K$, $t$ minimal and $\lambda_1=1$. Applying $y$ we get $\sum_i\sigma(\lambda_i)w_i=0$, and subtracting the two relations yields a shorter relation, which by the minimality of $t$ must be trivial. So $\sigma(\lambda_i)=\lambda_i$ for all $i$, that is, all $\lambda_i$ lie in $K_0$, which is a contradiction. Thus $\dim_{K_0}C_W(y)\le m$, and hence $|C_W(y)|\le|K_0|^m=Q^{m/s}$.
\end{proof}

The following elementary arithmetic fact will be useful.

\begin{lemma}\label{arith}
Let $p$ be a prime, $Q_0$ a prime power and $j\ge 0$.
\begin{enumerate}[(a)]
\item If $p$ divides $Q_0^{p^j}-1$, then $p$ divides $Q_0-1$.
\item If $p$ divides $Q_0-1$, then $p^j$ divides $(Q_0^{p^j}-1)/(Q_0-1)$.
\end{enumerate}
\end{lemma}

\begin{proof}
(a) The multiplicative order of $Q_0$ modulo $p$ divides both $p-1$ and $p^j$, so it is $1$.
(b) We have
$$\frac{Q_0^{p^j}-1}{Q_0-1}=\prod_{k=1}^j\left(1+R_k+R_k^2+\cdots+R_k^{p-1}\right),\qquad R_k=Q_0^{p^{k-1}}.$$
As $R_k\equiv 1\pmod p$, each factor is congruent to $p\equiv 0$ modulo $p$.
\end{proof}

Finally, we recall some standard facts on systems of imprimitivity. If $V$ is an irreducible $G$-module and $V=V_1\oplus\cdots\oplus V_t$ with $t>1$ is a decomposition into subspaces which are permuted transitively by $G$, we call it a system of imprimitivity of $V$, and we call the stabilizer $I$ of $V_1$ in $G$ the stabilizer of the system (it is determined up to conjugation).

\begin{lemma}\label{imprim}
Let $V$ be an irreducible $G$-module and $V=V_1\oplus\cdots\oplus V_t$ a system of imprimitivity with stabilizer $I$. Then the following hold.
\begin{enumerate}[(a)]
\item $V_1$ is an irreducible $I$-module and $V\cong\Ind_I^G(V_1)$.
\item If $Z\le G$ with $G=IZ$, then $V_Z\cong\Ind_{Z\cap I}^Z\big((V_1)_{Z\cap I}\big)$.
\end{enumerate}
\end{lemma}

\begin{proof}
(a) If $0<W<V_1$ were an $I$-submodule, then the sum of the subspaces $gW$, where $g$ runs over a transversal of $I$ in $G$, would be a proper nonzero $G$-submodule of $V$. The second statement is well known.
(b) As $G=IZ$, the group $Z$ permutes the $V_i$ transitively, and the stabilizer of $V_1$ in $Z$ is $Z\cap I$. Thus $V_Z$ is induced from $(V_1)_{Z\cap I}$.
\end{proof}

\section{Proof of Theorem B}

In this section, we prove Theorem B.

\begin{proof}[Proof of Theorem B]
Put $F=\GF(q^p)$ and let $N:F^\times\to\GF(q)^\times$ be the norm map, so that $U=\ker N$ is cyclic of order $1+q+\cdots+q^{p-1}$. We first show that $|U|_p=p$. If $p=2$, then $|U|=q+1\equiv 2\pmod 4$ by hypothesis. If $p$ is odd, write $q=1+mp$. Then
$$1+q+\cdots+q^{p-1}\equiv\sum_{i=0}^{p-1}(1+imp)=p+mp\cdot\frac{p(p-1)}{2}\equiv p\pmod{p^2}.$$
In both cases $|U|_p=p$. Since every element of $\GF(q)^\times$ of order $p$ has norm $\omega^p=1$, we have $U=\langle\omega\rangle\times K$, where $K$ is the Hall $p'$-subgroup of $U$. Moreover, $|U|\ge 1+q>p$, so that $K>1$.

The group $H=F\rtimes K$ is clearly solvable, and $|H|=q^p|K|$ is coprime to $p$. Write $t_{k,w}$ for the map $y\mapsto ky+w$. If $\pi\in P$ is given by $\pi(y)=\omega^iy^{q^j}$, then a direct computation shows that
$$\pi\, t_{k,w}\,\pi^{-1}=t_{k^{q^j},\pi(w)},$$
and since $K$ is characteristic in $F^\times$, this shows that $P$ acts on $H$ by conjugation. If $\pi$ acts trivially on $H$, then $\pi(w)=w$ for all $w\in F$, and taking $w=1$ we get $i\equiv 0\pmod p$. Since $F\neq\GF(q)$, we then also get $j\equiv 0\pmod p$. Thus $P$ acts faithfully on $H$.

We next determine the classes of $H$. Since $t_{k',u}\,t_{1,w}\,t_{k',u}^{-1}=t_{1,k'w}$, the class of a translation $t_{1,w}$ is $\{t_{1,kw}\mid k\in K\}$. If $k\neq 1$, then conjugating $t_{k,w}$ by the translations $t_{1,u}$ yields all elements $t_{k,w+(1-k)u}$, $u\in F$, and so the class of $t_{k,w}$ is $\{t_{k,w'}\mid w'\in F\}$.

Now let $c\in\cl(H)$. If $c=\{1\}$, then $C_P(c)=P$. If $c=\{t_{k,w'}\mid w'\in F\}$ with $k\neq 1$, then $\omega t_{k,w}\omega^{-1}=t_{k,\omega w}$ shows that $\omega\in C_P(c)$. Finally, let $c$ be the class of $t_{1,v}$ with $v\neq 0$. As $N(v)\in\GF(q)^\times$, we have $N(v^{q-1})=N(v)^{q-1}=1$, so $v^{q-1}\in U$, and we may write $v^{q-1}=\omega^{-i}k$ with $k\in K$. Put $\pi=\omega^i\varphi$. Then
$$\pi(v)=\omega^iv^q=\omega^iv^{q-1}v=kv,$$
and hence $\pi\,t_{1,v}\,\pi^{-1}=t_{1,kv}\in c$. Since $\pi\neq 1$, we have $|C_P(c)|\ge p$ also in this case. As $|P|=p^2$, the proof is complete.
\end{proof}

\begin{remark}\label{sharpC}
In the situation of Theorem B, assume that $q$ is a prime, and put $G_0=U\rtimes\langle\varphi\rangle$, acting on $V_0=F$ as a group of semilinear maps. Then $G_0$ is solvable and $p$-nilpotent with normal $p$-complement $K$, and $\langle\omega\rangle\le Z(G_0)$, so that $O_p(G_0)>1$. Since $K\cap\GF(q)^\times=1$ and $\GF(q)$ is the only proper subfield of $F$, the group $K$ generates $F$ as a field over the prime field $\GF(q)$, and so $V_0$ is a faithful irreducible $G_0$-module. Thus $(G_0,V_0)\in\X_p$. The computation in the proof of Theorem B shows that every $0\neq v\in V_0$ is fixed by the element $k^{-1}\omega^i\varphi\in G_0-U$ for suitable $k\in K$ and $i$, while $C_{G_0}(v)\cap U=1$, so that $C_{G_0}(v)$ embeds into $G_0/U\cong C_p$. Hence $|C_{G_0}(v)|_p=p=|G_0|_p^{1/2}$ for all $0\neq v\in V_0$. In particular, the inequality in Theorem C(a) is sharp.
\end{remark}

\section{A reduction}

In this section, we show that Theorem A follows from part (a) of Theorem C. 

\begin{lemma}\label{reduction}
Assume that $\mu_p(G,V)\le|G|_p^{1/2}$ for all $(G,V)\in\X_p$. Let $G$ be a solvable $p$-nilpotent group such that a Sylow $p$-subgroup $P$ of $G$ acts faithfully on $H=O_{p'}(G)$. Then there exists $x\in H$ such that $|C_G(x)|_p\le|G|_p^{1/2}$.
\end{lemma}

\begin{proof}
We argue by induction on $|G|$. If $P=1$, there is nothing to prove, so let $P>1$.

Step 1. $F(G)$ is a $p'$-group.

Since $O_p(G)$ and $H$ are normal subgroups of coprime orders, $[O_p(G),H]=1$. Thus $O_p(G)\le P$ acts trivially on $H$, and so $O_p(G)=1$ by hypothesis.

Step 2. We may assume that $\Phi(G)=1$.

Let $\Phi=\Phi(G)$ and suppose that $\Phi>1$. By Step 1, $\Phi\le F(G)\le H$. We show that $G/\Phi$ satisfies the hypotheses of the lemma. Clearly $O_{p'}(G/\Phi)=H/\Phi$. If $y\in P$ acts trivially on $H/\Phi$, then $y\Phi$ centralizes $F(G)/\Phi=F(G/\Phi)$ (see \cite[Theorem 1.12]{MW}), and as $G/\Phi$ is solvable, $y\Phi\in F(G/\Phi)$, that is, $y\in F(G)$. By Step 1, $y=1$. So by induction there is $x\in H$ with $|C_{G/\Phi}(x\Phi)|_p\le|G/\Phi|_p^{1/2}=|G|_p^{1/2}$. Since $C_G(x)\Phi/\Phi\le C_{G/\Phi}(x\Phi)$ and $\Phi$ is a $p'$-group, we have $|C_G(x)|_p\le|C_{G/\Phi}(x\Phi)|_p$, and we are done.

Step 3. We may assume that $p$ divides $|F(G/F(G))|$.

Let $F_0=1<F_1=F(G)<F_2<\cdots$ be the ascending Fitting series of $G$, and let $j\ge1$ be minimal such that $p$ divides $|F_{j+1}:F_j|$; such $j$ exists as $P>1$. Then $F_j$ is a $p'$-group. Let $S\le F_j$ be the normal subgroup of $G$ with $S/F_{j-1}=\Phi(F_j/F_{j-1})$. Suppose that $S>1$. Then $S\le H$ and $O_{p'}(G/S)=H/S$. If $y\in P$ acts trivially on $H/S$, then in $\bar G=G/F_{j-1}$ the element $\bar y$ centralizes $F(\bar G)/\Phi(F(\bar G))$. As $\Phi(F(\bar G))\le\Phi(\bar G)$, the element $\bar y$ also centralizes $F(\bar G)/\Phi(\bar G)$, and as in Step 2 we conclude that $\bar y\in F(\bar G)$, that is, $y\in F_j$, so $y=1$. Hence $G/S$ satisfies the hypotheses of the lemma, and as in Step 2 we are done by induction. So we may assume that $S=1$. Then $F_{j-1}=1$, so that $j=1$, and $p$ divides $|F_2:F_1|=|F(G/F(G))|$.

Step 4. The structure of $G$.

By Steps 2 and 3 and a theorem of Gasch\"utz (see \cite[Theorem 1.12]{MW}), $F:=F(G)$ is abelian, it is a direct product of minimal normal subgroups of $G$, and it has a complement $L$ in $G$. Since $G$ is solvable, $C_G(F)=F$, so that $L$ acts faithfully on $F$. Since $F$ acts trivially on itself by conjugation, the minimal normal subgroups of $G$ contained in $F$ are irreducible $L$-modules, so that $F$ is a faithful completely reducible $L$-module (possibly of mixed characteristic). Moreover $L\cong G/F$ is $p$-nilpotent and $O_p(L)>1$ by Step 3.

Step 5. Final conclusion.

Since $O_p(L)>1$ acts faithfully on $F$, there is a minimal normal subgroup $V_1$ of $G$ contained in $F$ on which $O_p(L)$ acts nontrivially. Let $V_2$ be the product of the remaining minimal normal subgroups in a decomposition of $F$, so that $F=V_1\times V_2$ and $V_2\trianglelefteq G$. Put $C=C_L(V_1)$ and $\bar L=L/C$. Then $(\bar L,V_1)\in\X_p$, because the image of $O_p(L)$ in $\bar L$ is a nontrivial normal $p$-subgroup of $\bar L$, and $V_1$ is a faithful irreducible $\bar L$-module over a prime field of characteristic different from $p$. Hence by hypothesis there exists $v_1\in V_1$ with
\begin{equation}\label{eq1}
|C_{\bar L}(v_1)|_p\le|\bar L|_p^{1/2}.
\end{equation}

If $C=1$, put $x=v_1$. As $F$ is abelian, $C_G(x)=F\,C_L(x)$, and thus $|C_G(x)|_p=|C_L(v_1)|_p\le|L|_p^{1/2}=|G|_p^{1/2}$, as desired. So we may assume that $C>1$.

Put $R=V_2C$. Since $C$ centralizes $V_1$ and normalizes $V_2$, it is easy to check that $R\trianglelefteq G$. Moreover, $R\cap V_1=1$, so that $R<G$. The group $R$ is solvable and $p$-nilpotent, and since $C$ acts faithfully on $V_2$ (as $L$ acts faithfully on $V_1\times V_2$ and $C$ centralizes $V_1$), a Sylow $p$-subgroup of $R$ acts faithfully on $V_2\le O_{p'}(R)$ (note that, as $V_2$ is abelian, every element $wc$ with $w\in V_2$ and $c\in C$ acts on $V_2$ in the same way as $c$). Thus by induction there exists $y\in O_{p'}(R)$ such that
\begin{equation}\label{eq2}
|C_R(y)|_p\le|R|_p^{1/2}.
\end{equation}
Note that $O_{p'}(R)$ is characteristic in $R$, hence normal in $G$, and so $O_{p'}(R)\le H$.

Now put $x=v_1y\in H$. If $g\in C_G(x)$, then $v_1y=v_1^gy^g$ yields $(v_1^g)^{-1}v_1=y^gy^{-1}\in V_1\cap O_{p'}(R)\le V_1\cap R=1$. Hence $C_G(x)=C_G(v_1)\cap C_G(y)$. Put $X=C_G(x)$. Then $|X|_p=|XR/R|_p\,|X\cap R|_p$, and $|X\cap R|_p\le|C_R(y)|_p\le|R|_p^{1/2}$ by \eqref{eq2}. On the other hand, $G/R$ is isomorphic to the semidirect product $V_1\rtimes\bar L$, and $XR/R$ centralizes the image of $v_1$. The centralizer of $v_1$ in $V_1\rtimes\bar L$ is $V_1\rtimes C_{\bar L}(v_1)$, and so $|XR/R|_p\le|C_{\bar L}(v_1)|_p\le|\bar L|_p^{1/2}$ by \eqref{eq1}. Altogether we obtain
$$|C_G(x)|_p\le\left(|\bar L|_p\,|R|_p\right)^{1/2}=\left(|L/C|_p\,|C|_p\right)^{1/2}=|L|_p^{1/2}=|G|_p^{1/2},$$
and the proof is complete.
\end{proof}

\begin{proof}[Proof of Theorem A, assuming Theorem C]
Let $G=HP$ be the semidirect product. Then $G$ is solvable and $p$-nilpotent with $O_{p'}(G)=H$, and $P$ acts faithfully on $H$. By Theorem C(a) with $N=G$, the hypothesis of Lemma \ref{reduction} is satisfied. Hence there exists $x\in H$ with $|C_G(x)|_p\le|G|_p^{1/2}=|P|^{1/2}$, and by Lemma \ref{classes} the class $c=x^H$ satisfies $|C_P(c)|=|C_G(x)|_p\le|P|^{1/2}$.
\end{proof}

\section{The imprimitive case}

In this section and the next one, we prove Theorem C by induction on $|G|$. Throughout, we let $(G,V)\in\X_p$, and we put $K=O_{p'}(G)$ and $Z=O_p(G)>1$. Note that $[K,Z]=1$, and that $Z$ is contained in every Sylow $p$-subgroup of $G$. We will use the following inductive hypothesis.

\begin{hypothesis}\label{IH}
Theorem C holds for all pairs $(\tilde G,\tilde V)\in\X_p$ with $|\tilde G|<|G|$ and all normal subgroups of $\tilde G$.
\end{hypothesis}

We will show that if $V$ is not quasiprimitive, then one of the following two configurations occurs.

\begin{definition}\label{configs}
Let $(G,V)\in\X_p$.
\begin{enumerate}[(I)]
\item We say that $(G,V)$ is in Configuration (I) if there is a system of imprimitivity $V=U_1\oplus\cdots\oplus U_p$ whose stabilizer $J$ is a normal subgroup of index $p$ in $G$, such that $O_p(J/C_J(U_1))>1$.
\item We say that $(G,V)$ is in Configuration (II) if there is a system of imprimitivity $V=V_1\oplus\cdots\oplus V_t$ whose stabilizer $I$ contains a Sylow $p$-subgroup $P$ of $G$ with $C_P(V_1)=1$.
\end{enumerate}
\end{definition}

Note that in Configuration (I), the pair $(J/C_J(U_1),U_1)$ lies in $\X_p$ by Lemma \ref{imprim}(a). In Configuration (II), we have $Z\le P\le I$ and $C_Z(V_1)\le C_P(V_1)=1$, so that the image of $Z$ in $I/C_I(V_1)$ is a nontrivial normal $p$-subgroup, and hence $(I/C_I(V_1),V_1)\in\X_p$ as well.

\begin{lemma}\label{coarsen}
Let $(G,V)\in\X_p$, and let $V=W_1\oplus\cdots\oplus W_s$ with $s>1$ be a system of imprimitivity such that the kernel $D$ of the permutation action of $G$ on $\{W_1,\dots,W_s\}$ contains $K$. Suppose that $Z\le D$ or that $V_Z$ is homogeneous. Then $(G,V)$ is in Configuration (I).
\end{lemma}

\begin{proof}
Since $K\le D$, the group $G/D$ is a $p$-group acting transitively on $\{W_1,\dots,W_s\}$. Let $\bar J$ be a maximal subgroup of $G/D$ containing the stabilizer of $W_1$; then $\bar J$ is normal of index $p$ in $G/D$. Let $J$ be its inverse image in $G$, and let $U_1$ be the sum of the $W_i$ in the $J$-orbit of $W_1$. Since $J$ is normal of index $p$ and contains the stabilizer of $W_1$, there are exactly $p$ such sums $U_1,\dots,U_p$, they are permuted regularly by $G/J$, and $V=U_1\oplus\cdots\oplus U_p$ is a system of imprimitivity with stabilizer $J$.

It remains to show that $O_p(J/C_J(U_1))>1$, and for this it suffices to show that $Z\cap J$ acts nontrivially on $U_1$. If $Z\le J$, this holds since otherwise the normal subgroup $Z$ would act trivially on every $gU_1$, $g\in G$, and hence on $V$. If $Z\not\le J$, then by hypothesis $V_Z$ is homogeneous, and as $G=JZ$, Lemma \ref{imprim}(b) shows that $V_Z\cong\Ind_{Z\cap J}^Z((U_1)_{Z\cap J})$. By Lemma \ref{permmod}, $Z\cap J$ cannot act trivially on $U_1$.
\end{proof}

\begin{proposition}\label{twoconfigs}
Let $(G,V)\in\X_p$ such that $V$ is not quasiprimitive. Then $(G,V)$ is in Configuration (I) or in Configuration (II).
\end{proposition}

\begin{proof}
Step 1. We may assume that $V_Z$ is homogeneous.

Otherwise, let $W_1,\dots,W_s$ be the homogeneous components of $V_Z$, so that $s>1$. As $Z\trianglelefteq G$, by Clifford's theorem these form a system of imprimitivity. Since $K$ centralizes $Z$, for $k\in K$ the subspace $kW_i$ is a $Z$-submodule isomorphic to $W_i$, and so $kW_i=W_i$. Thus $K$ and $Z$ lie in the kernel of the action of $G$ on $\{W_1,\dots,W_s\}$, and $(G,V)$ is in Configuration (I) by Lemma \ref{coarsen}.

Step 2. We may assume that $V_N$ is homogeneous for every $N\trianglelefteq G$ with $K\le N$.

Otherwise, the homogeneous components of $V_N$ form a system of imprimitivity whose kernel contains $N\ge K$, and as $V_Z$ is homogeneous by Step 1, Lemma \ref{coarsen} shows that $(G,V)$ is in Configuration (I).

Step 3. Final conclusion.

Since $V$ is not quasiprimitive, there exists $N\trianglelefteq G$ such that $V_N$ is not homogeneous. Let $V_1,\dots,V_t$ ($t>1$) be the homogeneous components of $V_N$, and let $I$ be the stabilizer of $V_1$. Since $K\trianglelefteq G$, the $K$-orbits on $\{V_1,\dots,V_t\}$ form a $G$-invariant partition. If $K$ is not transitive on $\{V_1,\dots,V_t\}$, then the sums of the $V_i$ over the $K$-orbits form a system of imprimitivity whose kernel contains $K$, and we are in Configuration (I) by Lemma \ref{coarsen}. So we may assume that $K$ is transitive on $\{V_1,\dots,V_t\}$. Then $t$ divides $|K|$, and hence $t=|G:I|$ is not divisible by $p$. Thus $I$ contains a Sylow $p$-subgroup $P$ of $G$.

Put $Q=C_P(V_1)$. Since $C_I(V_1)\trianglelefteq I$, we have $Q=P\cap C_I(V_1)\trianglelefteq P$. As $G=KP$, it follows that $KQ\trianglelefteq G$, and by Step 2, $V_{KQ}$ is homogeneous, say $V_{KQ}\cong eY$ for an irreducible $KQ$-module $Y$. By Clifford's theorem applied to $Q\trianglelefteq KQ$, the isomorphism types of the irreducible constituents of $Y_Q$ form a single orbit under conjugation by $KQ$. Since $Q$ acts trivially on $V_1\neq 0$, the trivial $Q$-module is a constituent of $V_Q\cong eY_Q$, hence of $Y_Q$. The trivial module is invariant under conjugation, and thus $Q$ acts trivially on $Y$, and hence on $V$. As $V$ is faithful, $Q=1$, and so $(G,V)$ is in Configuration (II).
\end{proof}

We now show that Theorem C holds in both configurations. We start with Configuration (II), which is easier.

\begin{theorem}\label{configII}
Assume Hypothesis \ref{IH}, and suppose that $(G,V)$ is in Configuration (II). Then Theorem C holds for $(G,V)$.
\end{theorem}

\begin{proof}
Let $V=V_1\oplus\cdots\oplus V_t$, $I$ and $P$ be as in Configuration (II), and put $\bar I=I/C_I(V_1)$. As observed above, $(\bar I,V_1)\in\X_p$, and $|\bar I|<|G|$ since $t>1$. Let $N\trianglelefteq G$, and let $\bar N$ be the image of $N\cap I$ in $\bar I$. Then $\bar N\trianglelefteq\bar I$.

If $0\neq w\in V_1$ and $n\in N$ with $nw=w$, then $nV_1$ is a component containing $w$, so $nV_1=V_1$. Thus $C_N(w)=C_{N\cap I}(w)$. Moreover, $C:=C_{N\cap I}(V_1)$ is a normal subgroup of $I$, so that $C\cap P$ is a Sylow $p$-subgroup of $C$. As $C\cap P\le C_P(V_1)=1$, the group $C$ is a $p'$-group. Since $C\le C_{N\cap I}(w)$, we conclude that
$$|C_N(w)|_p=|C_{\bar N}(w)|_p\quad\text{for all }0\neq w\in V_1.$$
Also, $N\cap P\le I$ is a Sylow $p$-subgroup of $N$, and so $|\bar N|_p=|N\cap I|_p=|N|_p$.

By Hypothesis \ref{IH} there is $w\in V_1$ such that $|C_{\bar N}(w)|_p\le|\bar N|_p^{1/2}=|N|_p^{1/2}$. If $|N|_p=1$, there is nothing to prove, and otherwise $w\neq 0$, and so $|C_N(w)|_p\le|N|_p^{1/2}$. This proves (a).

For (b), suppose that the action of $N$ on $V$ is tight, so that $|N|_p=p^{2a}$ with $a\ge1$ and $\mu_p(N,V)=p^a$. Since $|C_N(0)|_p=p^{2a}>p^a$, the computation above shows that the action of $\bar N$ on $V_1$ is tight as well, with the same minimum $p^a$. By Hypothesis \ref{IH}, there are minimizers $w,w'\in V_1$ of $\bar N$ which lie in different $\bar N$-orbits. Both are minimizers of $N$ on $V$. If $nw=w'$ for some $n\in N$, then $n\in N\cap I$ as above, and so $w'\in\bar Nw$, a contradiction. Thus the minimizers of $N$ lie in at least two $N$-orbits, which proves (b).
\end{proof}

We now turn to Configuration (I), which requires more work. For the rest of this section, let $(G,V)$ be in Configuration (I), and let $J$ and $V=U_1\oplus\cdots\oplus U_p$ be as in Definition \ref{configs}. Put $\bar J=J/C_J(U_1)$. Then $(\bar J,U_1)\in\X_p$ and $|\bar J|<|G|$. Let $N\trianglelefteq G$. We fix an element $x\in G-J$, and if $N\not\le J$, we choose $x\in N-J$. After renumbering, we may assume that $U_j=x^{j-1}U_1$ for $j=1,\dots,p$. Note that $x^p\in J$.

Put $N_J=N\cap J\trianglelefteq G$, and let $M_0=N_J$ and $M_j=C_{N_J}(U_1\oplus\cdots\oplus U_j)$ for $j=1,\dots,p$. Then $M_p=1$, all $M_j$ are normal in $J$, and $M_{j-1}/M_j$ acts faithfully on $U_j$. For $j=1,\dots,p$, let $\bar N_j$ be the image of $x^{-(j-1)}M_{j-1}x^{j-1}\le J$ in $\bar J$. Then $\bar N_j\trianglelefteq\bar J$, and the bijection $U_1\to U_j$, $w\mapsto x^{j-1}w$, transforms the action of $\bar N_j$ on $U_1$ into the action of $M_{j-1}/M_j$ on $U_j$. In particular,
\begin{equation}\label{layer}
|C_{M_{j-1}/M_j}(x^{j-1}w)|_p=|C_{\bar N_j}(w)|_p\quad\text{for } w\in U_1.
\end{equation}
Put $m_j=|\bar N_j|_p$, so that $m_1m_2\cdots m_p=|N_J|_p$. Since $x^{-(j-1)}M_{j-1}x^{j-1}\le x^{-(j-1)}N_Jx^{j-1}=N_J$, we have $\bar N_j\le\bar N_1$, and so $m_j\le m_1$ for all $j$.

For $w_1,\dots,w_p\in U_1$, we put
$$v(w_1,\dots,w_p)=\sum_{j=1}^px^{j-1}w_j\in V.$$
By Lemma \ref{chain} and \eqref{layer}, we have
\begin{equation}\label{chainbound}
|C_{N_J}(v(w_1,\dots,w_p))|_p\le\prod_{j=1}^p|C_{\bar N_j}(w_j)|_p.
\end{equation}

The $\bar N_1$-orbit of an element $w\in U_1$ will be called the type of $w$.

\begin{lemma}\label{rotation}
Suppose that $N\not\le J$ and that $C_N(v)\not\le N_J$, where $v=v(w_1,\dots,w_p)$. Then $w_1,\dots,w_p$ all have the same type.
\end{lemma}

\begin{proof}
Let $y\in C_N(v)-J$. Then $y=gx^i$ with $1\le i\le p-1$ and $g=yx^{-i}\in N\cap J=N_J$, since $x\in N$. Let $1\le j\le p$, and let $1\le j'\le p$ with $j'\equiv j+i\pmod p$. Since $y$ maps $U_j$ onto $U_{j'}$ and fixes $v$, we have $yx^{j-1}w_j=x^{j'-1}w_{j'}$. Put $n=j+i-1$. Then $n-(j'-1)\in\{0,p\}$, and we obtain
$$w_{j'}=x^{-(j'-1)}gx^{n}w_j=\left(x^{-(j'-1)}gx^{j'-1}\right)x^{n-(j'-1)}w_j.$$
Here $x^{-(j'-1)}gx^{j'-1}\in N_J$ as $N_J\trianglelefteq G$, and $x^{n-(j'-1)}\in\{1,x^p\}\subseteq N_J$. Hence $w_{j'}$ and $w_j$ have the same type. Since $i\not\equiv0\pmod p$, the map $j\mapsto j'$ is a $p$-cycle, and the assertion follows.
\end{proof}

\begin{theorem}\label{configI}
Assume Hypothesis \ref{IH}, and suppose that $(G,V)$ is in Configuration (I). Then Theorem C holds for $(G,V)$.
\end{theorem}

\begin{proof}
Let $N\trianglelefteq G$, and let notation be as above. For $j=1,\dots,p$, we call $w\in U_1$ $j$-good if $|C_{\bar N_j}(w)|_p\le m_j^{1/2}$, and we write $\G_j$ for the set of $j$-good vectors. By Hypothesis \ref{IH} applied to $(\bar J,U_1)$ and $\bar N_j$, every $\G_j$ is nonempty. Since $\bar N_j\trianglelefteq\bar J$, every $\G_j$ is a union of $\bar J$-orbits, and in particular a union of types. Also note that if $m_j=1$, then $\G_j=U_1$ contains at least two types, namely $\{0\}$ and a nonzero type.

Step 1. The case $N\le J$.

Then $N_J=N$. If we choose $w_j\in\G_j$ for all $j$, then by \eqref{chainbound} the vector $v=v(w_1,\dots,w_p)$ satisfies $|C_N(v)|_p\le\prod_jm_j^{1/2}=|N|_p^{1/2}$, which proves (a).

For (b), suppose that the action of $N$ on $V$ is tight, so that $\mu_p(N,V)=|N|_p^{1/2}$. For each $j$ let $c_j=\mu_p(\bar N_j,U_1)\le m_j^{1/2}$, and choose minimizers $w_j$ of $\bar N_j$ on $U_1$. By \eqref{chainbound}, $\mu_p(N,V)\le\prod_jc_j\le\prod_jm_j^{1/2}=\mu_p(N,V)$, so that $c_j=m_j^{1/2}$ for all $j$. In particular, the action of $\bar N_1$ on $U_1$ is either tight or $m_1=1$, and in both cases (using Hypothesis \ref{IH} in the first case) there are minimizers $w_1,w_1'$ of $\bar N_1$ on $U_1$ of different types. Put $v=v(w_1,w_2,\dots,w_p)$ and $v'=v(w_1',w_2,\dots,w_p)$. By \eqref{chainbound}, both $v$ and $v'$ are minimizers of $N$ on $V$. If $nv=v'$ for some $n\in N\le J$, then $n$ stabilizes $U_1$, and $nw_1=w_1'$, where $n$ acts on $U_1$ through its image in $\bar N_1$. This is a contradiction, and so the minimizers of $N$ lie in at least two $N$-orbits.

Step 2. The case $N\not\le J$: preliminaries.

Now $x\in N$, and since $N/N_J\cong NJ/J=G/J$, we have $|N|_p=p|N_J|_p$. Suppose first that we can choose $w_j\in\G_j$ for $j=1,\dots,p$ such that $w_1,\dots,w_p$ are not all of the same type. Then for $v=v(w_1,\dots,w_p)$ we have $C_N(v)\le N_J$ by Lemma \ref{rotation}, and by \eqref{chainbound},
$$|C_N(v)|_p\le\prod_jm_j^{1/2}=|N_J|_p^{1/2}=p^{-1/2}|N|_p^{1/2}<|N|_p^{1/2}.$$
Thus (a) holds, and the action of $N$ on $V$ is not tight, so that (b) holds as well.

So we may assume that no such choice is possible. Since every $\G_j$ is a nonempty union of types, this means that there is a single type $\Ob$ such that $\G_j=\Ob$ for $j=1,\dots,p$. In particular, $m_j>1$ for all $j$. For $w\in\Ob$ the number $c_j=|C_{\bar N_j}(w)|_p$ does not depend on the choice of $w$, as $\Ob$ is an $\bar N_1$-orbit and $\bar N_j\trianglelefteq\bar J$. Put
$$\delta=\prod_{j=1}^p\frac{m_j^{1/2}}{c_j}.$$
Each factor of $\delta$ is of the form $p^{k/2}$ with an integer $k\ge0$.

Step 3. We have $m_1^{1/2}/c_1>1$, and in particular $\delta\ge p^{1/2}$.

Suppose that $c_1=m_1^{1/2}$. By Hypothesis \ref{IH} there is $w\in U_1$ with $|C_{\bar N_1}(w)|_p\le m_1^{1/2}$, and every such $w$ lies in $\G_1=\Ob$ and thus satisfies $|C_{\bar N_1}(w)|_p=c_1$. Hence $\mu_p(\bar N_1,U_1)=m_1^{1/2}$, so that the action of $\bar N_1$ on $U_1$ is tight (recall $m_1>1$), and its minimizers form the single $\bar N_1$-orbit $\Ob$. Thus this action is rigid, contradicting Hypothesis \ref{IH}.

Step 4. Proof of (a) in the case $N\not\le J$.

Let $T$ be the set of all vectors $v(w_1,\dots,w_p)$ with $w_1,\dots,w_p\in\Ob$, and let $v\in T$. Since $|N:N_J|=p$, we have $|C_N(v)|_p\le p\,|C_{N_J}(v)|_p$, and so by \eqref{chainbound} and Step 3,
\begin{equation}\label{Tbound}
|C_N(v)|_p\le p\prod_jc_j=\frac{p}{\delta}\prod_jm_j^{1/2}=\frac{p^{1/2}}{\delta}|N|_p^{1/2}\le|N|_p^{1/2}.
\end{equation}
This proves (a).

Step 5. Proof of (b) in the case $N\not\le J$.

Suppose that the action of $N$ on $V$ is tight. By \eqref{Tbound} we then have $\delta\le p^{1/2}$, hence $\delta=p^{1/2}$ by Step 3, and by Step 3 again, the only factor of $\delta$ different from $1$ is the first one. Thus, writing $m_j=p^{e_j}$, we have that $e_1$ is odd and $c_1=p^{(e_1-1)/2}$. Moreover, by \eqref{Tbound} every $v\in T$ is a minimizer of $N$. If $v\in T$ satisfied $C_N(v)\le N_J$, then we would get $|C_N(v)|_p\le\prod_jc_j=p^{-1/2}|N|_p^{1/2}$ by \eqref{chainbound}, a contradiction. Hence $C_N(v)\not\le N_J$, and so $N=N_JC_N(v)$ and $Nv=N_Jv$ for all $v\in T$. Since $T$ is $N_J$-invariant (by the argument in the proof of Lemma \ref{rotation}), the $N$-orbits contained in $T$ coincide with the $N_J$-orbits contained in $T$. So if $N_J$ has at least two orbits on $T$, then (b) holds.

Thus suppose that $N_J$ is transitive on $T$, and fix $v\in T$. Then $|T|=|\Ob|^p=|N_J:C_{N_J}(v)|$. Write $|C_{N_J}(v)|_p=p^\gamma$. As $\Ob$ is an $\bar N_1$-orbit on which $|C_{\bar N_1}(w)|_p=p^{(e_1-1)/2}$, we have $|\Ob|_p=p^{(e_1+1)/2}$. Comparing $p$-parts, we get
$$e_1+\cdots+e_p-\gamma=\frac{p(e_1+1)}{2}.$$
On the other hand, $C_N(v)\not\le N_J$ implies $|C_N(v)|_p=p^{1+\gamma}$, and as $v$ is a minimizer of the tight action of $N$, we have $p^{1+\gamma}=|N|_p^{1/2}=p^{(1+e_1+\cdots+e_p)/2}$, that is, $\gamma=(e_1+\cdots+e_p-1)/2$. Combining the two equations yields
$$e_1+\cdots+e_p=pe_1+p-1>pe_1.$$
However, we have seen that $m_j\le m_1$, that is, $e_j\le e_1$ for all $j$, and so $e_1+\cdots+e_p\le pe_1$. This contradiction completes the proof.
\end{proof}

\section{The quasiprimitive case}

In this section, we prove Theorem C when $V$ is quasiprimitive. We first collect the facts on quasiprimitive solvable linear groups which we need; they are contained in \cite[Corollary 1.10, Corollary 2.3, Corollary 2.6]{MW}.

\begin{theorem}\label{structure}
Let $G$ be a solvable group and let $V$ be a faithful irreducible quasiprimitive $G$-module over $\GF(r)$. Let $F=F(G)$, let $Z$ be the socle of the cyclic group $Z(F)$, and put $A=C_G(Z)$. Then every normal abelian subgroup of $G$ is cyclic and acts fixed point freely on $V-\{0\}$, and there are normal subgroups $E$, $T$ and $U$ of $G$ with the following properties.
\begin{enumerate}[(a)]
\item $F=ET$, $E\cap T=Z=Z(E)$, $T=C_G(E)$ and $F=C_A(E/Z)$.
\item Every Sylow subgroup of $E$ is either cyclic of prime order or extraspecial, and the extraspecial Sylow subgroups of odd order have prime exponent.
\item $U$ is cyclic, $U\le T$, $|T:U|\le2$ and $C_T(U)=U$. In particular $Z\le Z(F)\le C_T(U)=U$.
\item $|E/Z|=e^2$ for an integer $e\ge1$, and $E/Z$ is a faithful $A/F$-module on which $A$ preserves the symplectic form given by commutation; in particular, if $E/Z$ has order $2^{2m}$, then $A/F$ embeds into $\Sp(2m,2)$.
\item If $E>Z$, then $E/Z=E_1/Z\times\cdots\times E_n/Z$ for chief factors $E_i/Z$ of $G$ with $Z(E_i)=Z$ and $|E_i/Z|=p_i^{2n_i}$ for primes $p_i$ and integers $n_i\ge1$.
\item If $e=1$, that is, if $E=Z$, then $G\le\Gamma(V)$.
\end{enumerate}
\end{theorem}

Here (e) is \cite[Corollary 1.10(iii),(iv)]{MW}, and (f) is \cite[Corollary 2.3(b)]{MW}, since $E=Z$ implies $F=T$. The invariant $e$ is the same as in \cite{HY}.

In the situation of Theorem \ref{structure}, the cyclic normal subgroup $U$ generates a subfield $\GF(Q)$ of $\mathrm{End}(V)$, where $Q=r^a$ for some $a\ge1$, and $|U|$ divides $Q-1$. Thus $V$ is a $\GF(Q)$-vector space, $B:=C_G(U)$ acts $\GF(Q)$-linearly on $V$, and as conjugation by elements of $G$ induces field automorphisms of $\GF(Q)$, the group $G/B$ embeds into the cyclic group $\Gal(\GF(Q)/\GF(r))$ of order $a$. Since $Z\le U$, we have $B\le A$. Moreover, $H=EU$ is a normal subgroup of $G$ with $Z(H)=U$ and $U\cap E=Z(E)$, and $V_H$ is homogeneous, say $V_H\cong tV_0$. By \cite[Corollary 2.6]{MW}, $\dim_{\GF(Q)}V_0=e$, so that $\dim_{\GF(r)}V=eat$. As the faithful irreducible characters of $H$ have degree $e$ (see the proof of \cite[Corollary 2.6]{MW}), the module $V_0$ is absolutely irreducible over $\GF(Q)$, that is, $\mathrm{End}_H(V_0)=\GF(Q)$, and so $\mathrm{End}_H(V)\cong M_t(\GF(Q))$. Finally note that $C_F(U)=E\,C_T(U)=EU$, since $E$ centralizes $T$, and so $B\cap F=EU$.

\begin{lemma}\label{splitting}
In the situation of Theorem \ref{structure}, the degree of $\mathrm{End}_{\GF(r)G}(V)$ over $\GF(r)$ divides $\dim_{\GF(r)}V/e$.
\end{lemma}

\begin{proof}
The field $K=\mathrm{End}_{\GF(r)G}(V)$ is contained in $\mathrm{End}_H(V)\cong M_t(\GF(Q))$, and it commutes with the center $\GF(Q)$ of this algebra. Hence $K$ and $\GF(Q)$ generate a subfield $L$ of $M_t(\GF(Q))$ containing $\GF(Q)$. As $\GF(Q)^t$ is a vector space over $L$, the degree $|L:\GF(Q)|$ divides $t$, and so $|K:\GF(r)|$ divides $|L:\GF(r)|$, which divides $at=\dim_{\GF(r)}V/e$.
\end{proof}

\subsection{The case $e=1$}

\begin{proposition}\label{semilin}
Let $(G,V)\in\X_p$ such that $V$ is quasiprimitive and $e=1$. Then we can identify $V$ with $\GF(|V|)$ in such a way that $G$ becomes a subgroup of the semilinear group $\Gamma(V)$, and $p$ divides $|V|-1$.
\end{proposition}

\begin{proof}
The first statement is Theorem \ref{structure}(f). Let $R$ be a minimal normal subgroup of $G$ contained in $O_p(G)>1$. Then $R$ is abelian, and so by Theorem \ref{structure} it is cyclic of order $p$ and acts fixed point freely on $V-\{0\}$. Thus $p$ divides $|V|-1$.
\end{proof}

The following result handles the semilinear case. Its proof uses Lemma \ref{arith}, and it is an elementary counting argument.

\begin{lemma}\label{gamma}
Let $Q$ be a prime power, let $p$ be a prime dividing $Q-1$, and let $N$ be a subgroup of $\Gamma(Q)=\{y\mapsto\alpha\sigma(y)\mid\alpha\in\GF(Q)^\times,\ \sigma\in\Gal(\GF(Q))\}$, acting on $V=\GF(Q)$. Let $C=N\cap\GF(Q)^\times$, and write $|C|_p=p^c$ and $|N:C|_p=p^b$. Then the following hold.
\begin{enumerate}[(a)]
\item $|C_N(v)|_p\le p^b$ for all $0\neq v\in V$.
\item There exists $0\neq v\in V$ with $|C_N(v)|_p\le p^{\min(b,c)}$. In particular, $\mu_p(N,V)\le|N|_p^{1/2}$.
\item The action of $N$ on $V$ is not rigid.
\end{enumerate}
\end{lemma}

\begin{proof}
(a) If $0\neq v\in V$, then $C_N(v)\cap C=1$, and so $C_N(v)$ is isomorphic to a subgroup of $N/C$, which embeds into the cyclic group $\Gal(\GF(Q))$. This proves (a), and shows that $C_N(v)$ is cyclic.

(b) By (a) we may assume that $c<b$. Let $D$ be the subgroup of order $p^{c+1}$ of the cyclic group $N/C$, and choose $g_0\in N$ such that $g_0C$ generates $D$. Write $g_0(y)=\alpha\sigma(y)$, so that $\sigma$ has order $p^{c+1}$, and let $\GF(Q_0)$ be the fixed field of $\sigma$, so that $Q=Q_0^{p^{c+1}}$.

Let $0\neq v\in V$ with $|C_N(v)|_p\ge p^{c+1}$. Since $C_N(v)$ embeds into the cyclic group $N/C$, its image contains $D$, so there is $h\in C_N(v)$ with $hC=g_0C$, that is, $h=zg_0$ with $z\in C$. Then $z\alpha\sigma(v)=v$, so $v/\sigma(v)\in\alpha C$.

Let $W$ be the kernel of the norm map from $\GF(Q)^\times$ to $\GF(Q_0)^\times$. The map $\theta(v)=v/\sigma(v)$ is a homomorphism from $\GF(Q)^\times$ onto $W$ by Hilbert's Theorem 90, and its kernel is $\GF(Q_0)^\times$. As $\alpha C\cap W$ is either empty or a coset of $C\cap W$, the number of $v$ as above is at most $|C\cap W|(Q_0-1)$, so that their proportion in $\GF(Q)^\times$ is at most
$$\frac{|C\cap W|}{|W|}\le\frac{|C\cap W|_p}{|W|_p}\le\frac{p^c}{|W|_p},$$
where we used that $W$ is cyclic. Now $|W|=(Q_0^{p^{c+1}}-1)/(Q_0-1)$, and since $p$ divides $Q-1=Q_0^{p^{c+1}}-1$, Lemma \ref{arith} shows that $p$ divides $Q_0-1$ and that $p^{c+1}$ divides $|W|$. Thus the proportion is at most $1/p<1$, and so there exists $0\ne v\in V$ with $|C_N(v)|_p\le p^c$. The last statement follows as $p^{\min(b,c)}\le p^{(b+c)/2}=|N|_p^{1/2}$.

(c) Suppose that the action of $N$ is tight. By (b), $p^{\min(b,c)}\ge\mu_p(N,V)=p^{(b+c)/2}$, so that $b=c\ge1$, and by (a) every nonzero vector is a minimizer. If the action were rigid, then $N$ would be transitive on $V-\{0\}$, so that $|N|=(Q-1)|C_N(v)|$ and hence $p^{2b}=|Q-1|_p\,p^b$, that is, $|Q-1|_p=p^b$. On the other hand, the image of $N$ in $\Gal(\GF(Q))$ has order divisible by $p^b$, so that $Q=Q_1^{p^b}$ for some prime power $Q_1$. As $p$ divides $Q-1$, Lemma \ref{arith} shows that $p$ divides $Q_1-1$ and that $p^{b+1}$ divides $Q-1$. This contradiction proves (c).
\end{proof}

\subsection{Groups with central $p'$-part}

\begin{lemma}\label{lemmaF}
Let $(G,V)\in\X_p$ such that $O_{p'}(G)\le Z(G)$. Then Theorem C holds for $(G,V)$. In particular, this is the case if $F(G)$ is a $p$-group.
\end{lemma}

\begin{proof}
Put $K=O_{p'}(G)$ and let $P$ be a Sylow $p$-subgroup of $G$. As $K$ is central and $G/K$ is a $p$-group, $P$ is normal in $G$, and $G=K\times P$. Let $N\trianglelefteq G$. Then $N_p=N\cap P$ is the normal Sylow $p$-subgroup of $N$, and hence $|C_N(v)|_p=|C_{N_p}(v)|$ for all $v\in V$. If $N_p=1$, there is nothing to prove. Otherwise $N_p$ acts faithfully and coprimely on $V$, and Theorem \ref{isaacs} yields $v\in V$ with $|C_{N_p}(v)|\le(|N_p|/p)^{1/p}<|N_p|^{1/2}$. Thus (a) holds, and the action of $N$ is not tight, so that (b) holds as well.

If $F(G)$ is a $p$-group, then $[K,F(G)]=1$, so that $K\le C_G(F(G))\le F(G)$ and hence $K=1$.
\end{proof}

The following observation is contained in the proof of \cite[Theorem 3.3]{SY} for $p=2$; the argument works for every prime.

\begin{lemma}\label{pnde}
Let $(G,V)\in\X_p$ such that $V$ is quasiprimitive, and use the notation of Theorem \ref{structure}. Then $p$ does not divide $e$.
\end{lemma}

\begin{proof}
Suppose that $p$ divides $e$. Then by Theorem \ref{structure}(e) there is a chief factor $E_i/Z$ of $G$ with $E_i\le E$ and $|E_i/Z|=p^{2n_i}$, where $n_i\ge1$. Put $K=O_{p'}(G)$. As $O_p(E)$ is a normal $p$-subgroup of $G$, we have $[K,O_p(E)]\le K\cap O_p(E)=1$, and since $E_i=(E_i\cap O_p(E))Z$, the group $K$ acts trivially on $E_i/Z$. Thus $G$ acts on the irreducible $\GF(p)G$-module $E_i/Z$ through the $p$-group $G/K$. But every irreducible module of a $p$-group over a field of characteristic $p$ is trivial, so that $|E_i/Z|=p$, a contradiction.
\end{proof}

\begin{lemma}\label{pnotdive}
Let $(G,V)\in\X_p$ such that $V$ is quasiprimitive, and use the notation of Theorem \ref{structure}. Suppose that $p$ is odd and that $p$ does not divide $|G:F(G)|$. Then $C_G(v)$ is a $p'$-group for every $0\neq v\in V$, and Theorem C holds for $(G,V)$.
\end{lemma}

\begin{proof}
As $p$ does not divide $|G:F(G)|$, a Sylow $p$-subgroup $P$ of $G$ is contained in $F(G)$, so that $P=O_p(G)$ is normal in $G$. As $p$ does not divide $e$ by Lemma \ref{pnde}, Theorem \ref{structure}(b) shows that the Sylow $p$-subgroup of $E$ is trivial or cyclic of order $p$, and since $|E:Z|=e^2$ is prime to $p$, it is contained in $Z\le T$. Hence $P\le T$, and as $p$ is odd and $|T:U|\le2$, we have $P\le U$. Since $U$ acts fixed point freely on $V-\{0\}$, for $0\neq v\in V$ the group $C_G(v)$ intersects the normal Sylow $p$-subgroup $P$ trivially, and thus is a $p'$-group. As observed after Definition \ref{tight}, this implies Theorem C for $(G,V)$.
\end{proof}

If $e\in\{2,3\}$, then $|G:F(G)|$ divides $|G:EU|=|G:B|\,|B:EU|$, where $|G:B|$ divides $a$ and $B/EU$ embeds into $A/F$, which embeds into $\Sp(2,e)$ of order $6$ or $24$. Hence in this case Lemma \ref{pnotdive} applies to every prime $p\ge5$ which does not divide $a$.

\subsection{The nonmetacyclic case}

Here we rely on the following result of D.\ F.\ Holt and Y.\ Yang \cite{HY}.

\begin{theorem}[Holt, Yang]\label{HY}
Let $G$ be a solvable group which acts faithfully, irreducibly and quasiprimitively on a finite vector space $V$. If $G$ is not metacyclic, then $G$ has a regular orbit on $V$, except for the cases listed in \cite[Table 4.1]{HY}.
\end{theorem}

In \cite[Table 4.1]{HY}, each case is given by the invariant $e$, the characteristic $r$ (denoted by $p$ there), $d=\dim_{\GF(r)}V$, an integer $a$ and $b=d/(ea)$, together with the number of groups without regular orbits and the largest order of such a group. If $e=1$, then $G$ is metacyclic by Theorem \ref{structure}(f), so that Theorem \ref{HY} only concerns the case $e>1$.

\begin{proposition}\label{nonmeta}
Let $(G,V)\in\X_p$ such that $V$ is quasiprimitive, $e>1$ and $G$ is not metacyclic. Then Theorem C holds for $(G,V)$.
\end{proposition}

\begin{proof}
If $G$ has a regular orbit on $V$, we are done, so by Theorem \ref{HY} we may assume that $(G,V)$ is one of the cases in \cite[Table 4.1]{HY}. In all of these cases, $e$ is a power of a prime $s$. We use the notation of Theorem \ref{structure} and the subsequent discussion. In particular, $|U|$ divides $r^a-1$, where $a$ divides $d/e$.

First consider the cases numbered 1, 3, 9, 10, 19, 20, 22, 25, 48, 62, 63, 65, 68 and 117 in \cite[Table 4.1]{HY}. A direct check shows that in each of them, $r^{a'}-1$ is a power of $s$ for every divisor $a'$ of $d/e$. Hence $U$ is an $s$-group. If $s=2$, then $T$ is a $2$-group, since $|T:U|\le2$. If $s$ is odd, then $r^{d/e}-1$ is odd, so that $r=2$ and $U$ has odd order; then $|T|_2\le|T:U|\le2$, and as $T\le F$ is nilpotent, $T$ is cyclic, so that $T=C_T(U)=U$. In both cases $T$ is an $s$-group. The Sylow subgroups of $E$ which are cyclic of prime order are contained in $Z\le T$, and so $F(G)=ET$ is an $s$-group. Since $1<O_p(G)\le F(G)$, we have $p=s$, which is impossible by Lemma \ref{pnde}. So these cases do not occur for pairs in $\X_p$.

In the remaining $17$ cases, listed in Table \ref{table1}, $e$ is again a power of a prime $s$. As $1<O_p(G)\le F(G)$, $|F(G)|=e^2|T|$, $|T:U|\le2$ and $p\nmid e$ by Lemma \ref{pnde}, the prime $p$ divides $2|U|$, where $|U|$ divides $r^{a'}-1$ for a divisor $a'$ of $d/e$. If $e\in\{2,3\}$ and $p\ge5$ does not divide $d/e$, then $p$ does not divide $a$, as $a$ divides $d/e$, and we are done by Lemma \ref{pnotdive} and the subsequent remark. The pairs consisting of a case and a prime $p$ which are not covered by this are listed in the last column of Table \ref{table1}. These pairs were checked by computer; in fact, our computations covered all primes $p$ in all $17$ cases. We used GAP 4.13.1 \cite{GAP} with the IRREDSOL library \cite{IRREDSOL} (see Appendix \ref{appendix} for the code and a summary of the results), which contains, up to conjugacy, all irreducible solvable subgroups of $\GL(d,r)$ for $r^d<2^{24}$. In all remaining cases we have $d=ea$ in \cite[Table 4.1]{HY}, and so by Lemma \ref{splitting} the degree of the splitting field of $V$ over $\GF(r)$ divides $d/e$. For each case, we considered all irreducible solvable subgroups $G$ of $\GL(d,r)$ whose splitting field has degree dividing $d/e$ and whose order does not exceed the largest order given in \cite[Table 4.1]{HY}; see Table \ref{table1} below. For each such $G$ and each prime $p\neq r$ such that $G$ is $p$-nilpotent and $O_p(G)>1$, we proceeded as follows. If $F(G)$ is a $p$-group, Lemma \ref{lemmaF} applies. Otherwise, we tested random vectors $v\in V$ until we found one for which $C_G(v)$ is a $p'$-group; as observed after Definition \ref{tight}, this proves Theorem C for $(G,V)$. If no such vector was found after $25$ attempts, we computed all orbits of $G$ on $V$ together with their stabilizers, and for every normal subgroup $N$ of $G$ we computed $\mu_p(N,V)$ and, in the tight case, the number of $N$-orbits of minimizers, which equals the sum of the indices $|G:NC_G(v)|$ over representatives $v$ of the $G$-orbits of minimizers. In all cases, Theorem C was confirmed.

In the largest case, number 28 with $e=4$ and $V=\GF(5)^8$, we additionally skipped all groups violating one of the following necessary conditions for the groups which are not covered by Lemma \ref{lemmaF}. Here $|U|$ divides $24$. If $3$ does not divide $|U|$, then as above $F(G)$ is a $2$-group and Lemma \ref{lemmaF} applies. Otherwise $Q=25$, and $6$ divides $|U|$ as $|Z(E)|$ is even. Then $|EU|=e^2|U|=16|U|$ is divisible by $96$ and divides $384$. As $B\cap F=EU$ and $B\le A$, the group $B/EU$ embeds into $A/F$, which embeds into $\Sp(4,2)$ by Theorem \ref{structure}(d), and $|G:B|\le2$. Thus $96$ divides $|G|$, and $|G|$ divides $2\cdot|\Sp(4,2)|\cdot384=552960$. If the splitting field of $V$ is $\GF(5)$, then $G\neq B$, as otherwise $\GF(25)\subseteq\mathrm{End}_{\GF(5)G}(V)$, and so $192$ divides $|G|$. Finally, $F(G)$ is nonabelian with cyclic center, and $32$ divides $|F(G)|$.
\end{proof}

\begin{table}[ht]
\centering
\begin{tabular}{ccccccc}
\toprule
Case in \cite[Table 4.1]{HY} & $e$ & $r$ & $d=\dim V$ & $a$ & largest $|G|$ & primes $p$ to check\\
\midrule
21 & 4 & 7 & 4 & 1 & 6912 & 3\\
23 & 4 & 11 & 4 & 1 & 11520 & 5\\
24 & 4 & 13 & 4 & 1 & 13824 & 3\\
28 & 4 & 5 & 8 & 2 & 55296 & 3\\
49 & 3 & 7 & 3 & 1 & 1296 & 2\\
50 & 3 & 13 & 3 & 1 & 2592 & 2\\
51 & 3 & 2 & 12 & 4 & 12960 & --\\
52 & 3 & 19 & 3 & 1 & 3888 & 2\\
53 & 3 & 5 & 6 & 2 & 10368 & 2\\
64 & 2 & 7 & 2 & 1 & 144 & 3\\
66 & 2 & 11 & 2 & 1 & 240 & --\\
67 & 2 & 13 & 2 & 1 & 288 & 3\\
69 & 2 & 19 & 2 & 1 & 432 & 3\\
71 & 2 & 5 & 4 & 2 & 1152 & 3\\
72 & 2 & 3 & 6 & 3 & 1872 & --\\
74 & 2 & 7 & 4 & 2 & 2304 & 3\\
75 & 2 & 3 & 8 & 4 & 7680 & --\\
\bottomrule
\end{tabular}
\caption{The cases of \cite[Table 4.1]{HY} treated by computer.}\label{table1}
\end{table}

\begin{remark}\label{exhaustive}
For small parameters, we actually checked all irreducible solvable subgroups of $\GL(d,r)$, without any restriction on the order or on the splitting field. More precisely, Theorem C was confirmed for all $(G,V)\in\X_p$ with $G$ an irreducible solvable subgroup of $\GL(d,r)$ for $d=2$ and $r\in\{3,5,7,11,13,17,19,23\}$, for $d=3$ and $r\in\{5,7,13,19\}$, for $d=4$ and $r\in\{3,5,7\}$, and for $d\in\{5,6\}$ and $r\in\{2,3\}$. This covers the cases 21, 49, 50, 52, 64, 66, 67, 69, 71, 72 and 74 of Table \ref{table1}, and these computations also serve as an independent check of the arguments in Section 5. In all computations, more than $1200$ pairs required a complete orbit analysis, and among those, more than $240$ had a normal subgroup with a tight action; in each of these, the minimizers formed at least two orbits.
\end{remark}

\subsection{The metacyclic case}

Now let $(G,V)\in\X_p$ such that $V$ is quasiprimitive, $e>1$ and $G$ is metacyclic. We use the notation of Theorem \ref{structure} and the subsequent discussion; in particular $B=C_G(U)$ and $\GF(Q)$ is the field generated by $U$, where $Q=r^a$.

\begin{lemma}\label{metastruct}
In this situation, $e=2$, the Sylow $2$-subgroup of $E$ is extraspecial of order $8$, $p$ and $r$ are odd, $|B:EU|\le2$, $|B|\le8(Q-1)$ and $|G|\le8a(Q-1)$. Moreover, $V$ is a $\GF(Q)$-space of dimension $2k$ for some $k\ge1$, and $B$ is an $r'$-group.
\end{lemma}

\begin{proof}
A metacyclic group $X$ of exponent $s$ has order at most $s^2$, and it is generated by two elements. By Theorem \ref{structure}(b), the extraspecial Sylow subgroups of $E$ of odd order have prime exponent, so by the first fact there are none. Since $e>1$, the Sylow $2$-subgroup $E_2$ of $E$ is extraspecial, and by the second fact $|E_2/\Phi(E_2)|\le4$. As $\Phi(E_2)=Z(E_2)$, we get $|E_2|=8$, so that $e=2$ and $E/Z\cong E_2/Z(E_2)$ has order $4$. As $B\cap F=EU$ and $B\le A$, the group $B/EU$ embeds into $A/F$, which embeds into $\Sp(2,2)\cong S_3$ by Theorem \ref{structure}(d). If $3$ divided $|B:EU|$, we could choose a $3$-element $y\in B$ acting nontrivially, and hence fixed point freely, on $E/Z$. Then $[E,y]Z=E$. If $[E,y]$ were cyclic, then $E=[E,y]Z$ would be abelian, as $Z=Z(E)$. Hence $[E,y]$ is not cyclic, and so the derived subgroup of the metacyclic group $E\langle y\rangle$ would not be cyclic. This is impossible, as metacyclic groups have cyclic derived subgroups. Thus $|B:EU|\le2$. Since the Sylow $2$-subgroup $E_2$ of $E$ is a nontrivial normal $2$-subgroup of $G$, and $V$ is faithful and irreducible, we have $r\neq2$, and $p\neq2$ by Lemma \ref{pnde}. As $|EU|=e^2|U|=4|U|$ and $|U|$ divides $Q-1$, we get $|B|\le8(Q-1)$, and $B$ is an $r'$-group. Also $|G:B|\le a$. Finally, $\dim_{\GF(r)}V=eat=2at$, so that $\dim_{\GF(Q)}V=2t$, and we put $k=t$.
\end{proof}

\begin{lemma}\label{metafix}
Let $x\in G$ be of prime order $s$ with $x\notin U$.
\begin{enumerate}[(a)]
\item If $x\in B$, then $|C_V(x)|\le Q^k$.
\item If $x\notin B$, then $x$ acts $\sigma$-semilinearly on $V$ for some $\sigma\in\Gal(\GF(Q))$ of order $s$, and $|C_V(x)|\le Q^{2k/s}$.
\end{enumerate}
\end{lemma}

\begin{proof}
(a) Recall that $V_H\cong kV_0$ for the absolutely irreducible $\GF(Q)H$-module $V_0$ of dimension $2$, where $H=EU$. Write $V=V_0\otimes M$ with $\dim_{\GF(Q)}M=k$. The $\GF(Q)$-span of $H$ in $\mathrm{End}(V)$ is $M_2(\GF(Q))\otimes1$. Since $x$ normalizes $H$ and centralizes $U$, conjugation by $x$ induces a $\GF(Q)$-algebra automorphism of $M_2(\GF(Q))$, which is inner by the Skolem--Noether theorem, so there is $w\in\GL(V_0)$ such that $(w\otimes1)^{-1}x$ centralizes $M_2(\GF(Q))\otimes1$. Hence $x=w\otimes m$ for some $m\in\GL(M)$. If $w$ were a scalar, then $x$ would centralize $E$, so that $x\in C_G(E)\cap B=C_T(U)=U$ by Theorem \ref{structure}, which is not the case. As $x$ has prime order $s\neq r$ (note that $B$ is an $r'$-group by Lemma \ref{metastruct}), we have $(w\otimes m)^s=1$, and hence $w^s=\lambda\cdot1$ for a nonzero scalar $\lambda$. Thus the minimal polynomial of $w$ divides $X^s-\lambda$, which has no repeated roots in characteristic $r$, and so $w$ is semisimple. Being semisimple and non-scalar, $w$ has two distinct eigenvalues $\alpha\neq\beta$ over the algebraic closure. Then the fixed space of $x$ has dimension equal to the sum of the multiplicities of $\alpha^{-1}$ and $\beta^{-1}$ as eigenvalues of $m$, which is at most $k$.

(b) As $x\notin B$, the coset $xB$ has order $s$ in $G/B$, and so $x$ is $\sigma$-semilinear for some $\sigma$ of order $s$. Now apply Lemma \ref{semilinear}.
\end{proof}

\begin{proposition}\label{metak2}
In the situation of this subsection, suppose that $k\ge2$. Then $G$ has a regular orbit on $V$, unless $(Q,k)\in\{(3,2),(5,2),(7,2),(9,2)\}$.
\end{proposition}

\begin{proof}
If $G$ has no regular orbit, then every vector is fixed by an element of prime order. Elements of $U$ fix only $0$, and every other element of prime order has at most $Q^k$ fixed points by Lemma \ref{metafix}. Since there are fewer than $|G|\le8a(Q-1)$ such elements, we get $Q^{2k}\le8a(Q-1)Q^k$, that is, $Q^k\le8a(Q-1)$. For $k\ge3$ this would give $Q^2<8a$, which is impossible for $Q=r^a$ with $r\ge3$. For $k=2$, we get $Q<8a$, and as $r\ge3$, this leaves only $Q\in\{3,5,7,9\}$.
\end{proof}

For the case $k=1$, we fix $N\trianglelefteq G$ and write $|N\cap B|_p=p^\gamma$ and $|NB:B|_p=p^\beta$, so that $|N|_p=p^{\beta+\gamma}$. Note that $NB/B$ is cyclic. The following counting argument is similar to the proof of Lemma \ref{gamma}(b).

\begin{lemma}\label{metacount}
Suppose that $k=1$, let $1\le j\le\beta$, and let $\Omega$ be a set of vectors $v\in V$ such that $C_{N\cap B}(v)$ is a $p'$-group. Then the number of $v\in\Omega$ with $|C_N(v)|_p\ge p^j$ is at most $8(Q-1)Q^{2/p^j}$.
\end{lemma}

\begin{proof}
Let $v\in\Omega$ with $|C_N(v)|_p\ge p^j$. Then a Sylow $p$-subgroup of $C_N(v)$ intersects $B$ trivially, and so it embeds into the cyclic group $NB/B$. Hence its image contains the unique subgroup $D$ of order $p^j$ of $NB/B$. Fix $g_0\in N$ with $\langle g_0B\rangle=D$. Then there is $h\in C_N(v)$ with $hB=g_0B$, that is, $h=zg_0$ with $z\in N\cap B$. By Lemma \ref{semilinear}, $h$ has at most $Q^{2/p^j}$ fixed points. As there are at most $|N\cap B|\le|B|\le8(Q-1)$ choices for $z$, the lemma follows.
\end{proof}

\begin{proposition}\label{metak1odd}
In the situation of this subsection, suppose that $k=1$ and that $p$ is odd. Then Theorem C holds for $(G,V)$, unless $a=3$ and $r\in\{5,7\}$.
\end{proposition}

\begin{proof}
Let $N\trianglelefteq G$. Since $|B:U|\le8$, the Sylow $p$-subgroup of $B$ is contained in $U$, which acts fixed point freely on $V-\{0\}$. So $C_{N\cap B}(v)$ is a $p'$-group for all $0\neq v\in V$, and as in the proof of Lemma \ref{metacount}, $|C_N(v)|_p\le p^\beta$ for all $0\neq v\in V$.

If $\beta=0$, then $\mu_p(N,V)=1$, and we are done. So let $\beta\ge1$; then $p$ divides $|G:B|$, which divides $a$, and hence $a\ge3$.

Suppose first that $\gamma\ge\beta$. Then $|C_N(v)|_p\le p^\beta\le|N|_p^{1/2}$ for all $0\neq v\in V$, which proves (a). If the action of $N$ is tight, then all nonzero vectors are minimizers. If it were rigid, then $N$ would be transitive on $V-\{0\}$, so that $Q^2-1\le|N|\le8(Q-1)a$, that is, $Q+1\le8a$. But $Q\ge3^a$ and $a\ge3$, so $Q+1>8a$, a contradiction.

Now suppose that $\gamma<\beta$, and put $j=\lfloor(\beta+\gamma)/2\rfloor+1$. Then $1\le j\le\beta$, and by Lemma \ref{metacount} (with $\Omega=V-\{0\}$), at most $8(Q-1)Q^{2/p}$ nonzero vectors $v$ satisfy $|C_N(v)|_p>p^{j-1}$. Every other nonzero vector satisfies $|C_N(v)|_p\le p^{\lfloor(\beta+\gamma)/2\rfloor}\le|N|_p^{1/2}$. Hence (a) holds if $Q^2-1>8(Q-1)Q^{2/p}$. Moreover, if the action of $N$ is tight, then all these vectors are minimizers, and they lie in at least two $N$-orbits as soon as their number exceeds $|N|\le8(Q-1)a$. Thus Theorem C holds for $N$ provided that
$$Q+1>8\left(Q^{2/p}+a\right).$$
Since $p$ divides $a$ and $r\neq p$ is odd, a direct check shows that this inequality fails only for $p=a=3$ and $r\in\{5,7\}$.
\end{proof}

\begin{proposition}\label{meta}
Let $(G,V)\in\X_p$ such that $V$ is quasiprimitive, $e>1$ and $G$ is metacyclic. Then Theorem C holds for $(G,V)$.
\end{proposition}

\begin{proof}
Since $p$ is odd by Lemma \ref{metastruct}, Propositions \ref{metak2} and \ref{metak1odd} show that it remains to consider finitely many exceptional cases. By Lemma \ref{metastruct}, $V=\GF(r)^d$ with $d=2ak$ and $|G|\le8a(Q-1)$, and by Lemma \ref{splitting} the degree of the splitting field of $V$ divides $d/2=ak$. The exceptional cases are:
\begin{enumerate}[(i)]
\item $k=2$ and $Q\in\{3,5,7\}$, that is, subgroups of $\GL(4,r)$ with $r\in\{3,5,7\}$;
\item $k=2$ and $Q=9$, that is, subgroups of $\GL(8,3)$ of order at most $128$;
\item $k=1$, $a=3$, $r\in\{5,7\}$, that is, subgroups of $\GL(6,5)$ and $\GL(6,7)$ of order at most $2976$ and $8208$, respectively.
\end{enumerate}
The groups in (i) are covered by Remark \ref{exhaustive}. The groups in (ii) are covered by the computations for case 75 in Table \ref{table1}, in which all irreducible solvable subgroups of $\GL(8,3)$ of order at most $7680$ whose splitting field has degree $1$, $2$ or $4$ were checked. For (iii), we checked all irreducible solvable subgroups of $\GL(6,5)$ of order at most $2976$ and of $\GL(6,7)$ of order at most $8208$ whose splitting field has degree $1$ or $3$, in the same way as in the proof of Proposition \ref{nonmeta}. In all cases, Theorem C was confirmed.
\end{proof}

\section{Proofs of Theorem C and Corollary D}

We are now ready to complete the proofs of our main results.

\begin{proof}[Proof of Theorem C]
We argue by induction on $|G|$, so that Hypothesis \ref{IH} holds. If $V$ is not quasiprimitive, then by Proposition \ref{twoconfigs}, $(G,V)$ is in Configuration (I) or (II), and we are done by Theorems \ref{configI} and \ref{configII}. So let $V$ be quasiprimitive, and let $e$ be as in Theorem \ref{structure}. If $e=1$, then by Proposition \ref{semilin}, $G\le\Gamma(V)$ and $p$ divides $|V|-1$, and Theorem C follows from Lemma \ref{gamma}, applied to every normal subgroup $N$ of $G$. If $e>1$, Theorem C follows from Propositions \ref{nonmeta} and \ref{meta}.
\end{proof}

Since Theorem A has already been derived from Theorem C in Section 4, it remains to prove Corollary D.

\begin{proof}[Proof of Corollary D]
For every subgroup $B\le P$, the number of $B$-invariant classes of $H$ equals the number of $B$-invariant irreducible characters of $H$, since the action is coprime (see \cite[Chapter 13]{IsCT}). By a well-known theorem of Burnside on marks (see, e.g., \cite[Proposition 1.2.2]{tD}), two finite $P$-sets are isomorphic if every subgroup of $P$ has the same number of fixed points on them. Hence the actions of $P$ on $\cl(H)$ and on $\Irr(H)$ are permutation isomorphic, and Corollary D follows from Theorem A.
\end{proof}

We conclude this section by showing that part (b) of Theorem C is not an artifact of our proof, but a consequence of part (a).

\begin{proposition}\label{wreath}
Let $(H,U)\in\X_p$ such that the action of $H$ on $U$ is rigid. Then $G=H\wr C_p$, acting naturally on $V=U^p$, satisfies $(G,V)\in\X_p$ and $\mu_p(G,V)>|G|_p^{1/2}$.
\end{proposition}

\begin{proof}
Clearly $G$ is solvable, it is $p$-nilpotent with normal $p$-complement $O_{p'}(H)^p$, and $O_p(H)^p\le O_p(G)$. Since the components $U_i$ of $V$ are irreducible modules for the base group $H^p$ with pairwise different kernels, and since they are permuted transitively, $V$ is a faithful irreducible $G$-module. Thus $(G,V)\in\X_p$.

Let $|H|_p=p^{2a}$ with $a\ge 1$, so that $|G|_p=p^{2ap+1}$ and $|G|_p^{1/2}=p^{ap+1/2}$, and let $\Ob$ be the $H$-orbit of minimizers on $U$. Let $v=(u_1,\dots,u_p)\in V$. If some $u_i\notin\Ob$, then $|C_H(u_i)|_p\ge p^{a+1}$, and as $C_G(v)\ge\prod_iC_H(u_i)$, we get $|C_G(v)|_p\ge p^{ap+1}$. If $u_i\in\Ob$ for all $i$, then $H^p$ is transitive on $\Ob^p$, so if $x$ is a generator of the top group $C_p$, there is $h\in H^p$ with $xv=hv$. Then $h^{-1}x\in C_G(v)-H^p$, and so $|C_G(v)|_p=p\,|C_{H^p}(v)|_p=p^{ap+1}$. In both cases $|C_G(v)|_p>|G|_p^{1/2}$.
\end{proof}

Thus Theorem C(a) for all pairs in $\X_p$ implies that no action in $\X_p$ is rigid, so that Theorem C(b) holds at least for $N=G$.

\section{Outlook}

We conclude this paper with some remarks and open questions.

First, it would be desirable to have a proof of Theorem C which does not depend on the classification in \cite{HY} and on computer calculations. By Lemma \ref{pnde}, only the case where $p$ does not divide $e$ has to be considered there. Here the elementary bound $|C_G(v)|_p\le|G|_p/|O_p(Z(G))|$ for $0\neq v\in V$, which follows from the fixed point free action of $O_p(Z(G))$, together with character-theoretic bounds for fixed points settles some cases, but generic counting arguments lose too much in exactly the small cases occurring in \cite[Table 4.1]{HY}, and eliminating the computations would require a new structural idea.

Second, it is natural to ask whether the solvability hypothesis on $H$ in Theorem A is necessary. Our reduction in Section 4 uses the Fitting series of $H$ in an essential way, and, as pointed out in the introduction, the Hartley--Turull lemma is not available for classes.

Third, Isaacs asked in \cite[Section 5]{Is99} whether for every coprime action of a group $N$ on a group $H$ there exists $x\in H$ with $|C_N(x)|\le|N|^{1/2}$. This was answered in the affirmative by Z.\ Halasi and K.\ Podoski \cite{HP}, who proved that every coprime linear group $G\le\GL(V)$ has a base of size two, that is, there are $u,v\in V$ with $C_G(u)\cap C_G(v)=1$; this gives $|G|\le|u^G|\,|v^G|$, and so one of these orbits has size at least $|G|^{1/2}$. This suggests the following question: If $N$ is a group acting faithfully and coprimely on a solvable group $H$, does there always exist a class $c$ of $H$ with $|C_N(c)|\le|N|^{1/2}$? Theorem A answers this question for $p$-groups $N$, and Theorem B shows that for classes a bound better than $|N|^{1/2}$ cannot be expected in general. By the same argument as above, a positive answer would follow if there always exist two classes $c,c'$ of $H$ with $C_N(c)\cap C_N(c')=1$. Note that in the examples of Theorem B such classes exist, for instance a class $\{t_{k,w'}\mid w'\in F\}$ with $k\neq1$, whose stabilizer is $\langle\omega\rangle$, together with a class of a translation, whose stabilizer is of the form $\langle\omega^i\varphi\rangle$.

Finally, the examples in Theorem B suggest the problem of classifying all coprime actions of $p$-groups $P$ on solvable groups $H$ for which $\min_{c\in\cl(H)}|C_P(c)|=|P|^{1/2}$, similar in spirit to the classification of the linear group actions with $|G:G'|$ equal to the largest orbit size in \cite{CK}.

\appendix

\section{GAP code and summary of the computations}\label{appendix}

The computations described in Section 6 were carried out with GAP 4.13.1 \cite{GAP} and the development version 1.4.dev of the IRREDSOL package \cite{IRREDSOL}, obtained from its GitHub repository. In this appendix we reproduce the code and summarize the results. The code below is the version used for the largest runs. For the exhaustive runs and for the runs concerning $\GL(4,11)$, $\GL(4,13)$ and the groups of $\GL(6,5)$ with splitting field $\GF(5)$, we used an earlier version of the main loop, which differs only in that it does not skip the pairs $(G,p)$ covered by Lemma \ref{lemmaF} and does not apply the additional filters for $\GL(8,5)$; thus it checks more pairs.

For a group $G$ from the library and a prime $p$, the code first tests whether $G$ is $p$-nilpotent with $O_p(G)>1$, so that $(G,V)\in\X_p$. It then looks for a vector $v$ such that $C_G(v)$ is a $p'$-group, which proves Theorem C for $(G,V)$ (see the remark after Definition \ref{tight}). If no such vector is found among $25$ random vectors, the function \texttt{FullAnalysis} computes all orbits of $G$ on $V-\{0\}$ and checks parts (a) and (b) of Theorem C for every normal subgroup $N$ of $G$ with $|N|_p>1$. Here the number of $N$-orbits of minimizers is computed as the sum of $|G:NC_G(v)|$ over representatives $v$ of the $G$-orbits of minimizers, as in the proof of Proposition \ref{nonmeta}. The vector $0$ never needs to be considered, since $|C_N(0)|_p=|N|_p$.

\begin{lstlisting}
LoadPackage("irredsol");;

LOG := "results.log";;
Out := function(arg)
  AppendTo(LOG, Concatenation(List(arg, String)), "\n");
end;;

# the l-part of the integer n
lpart := function(n, l)
  local r;
  r := 1;
  while n mod l = 0 do n := n/l; r := r*l; od;
  return r;
end;;

# H is a pc group acting on F^n, where pcgs is a pcgs of H and
# mats are the matrices by which the elements of pcgs act.
# Checks Theorem C (a) and (b) for all normal subgroups N of H
# with |N|_l > 1.
FullAnalysis := function(H, pcgs, mats, F, n, l, tag)
  local V, orbs, reps, stabs, N, nl, vals, m, mins, cnt,
        nviol, ntight, nrigid, nn;
  V := Filtered(Elements(F^n), v -> not IsZero(v));
  orbs := OrbitsDomain(H, V, pcgs, mats, OnRight);
  reps := List(orbs, o -> o[1]);
  stabs := List(reps, v -> Stabilizer(H, v, pcgs, mats, OnRight));
  nviol := 0; ntight := 0; nrigid := 0; nn := 0;
  for N in NormalSubgroups(H) do
    nl := lpart(Size(N), l);
    if nl = 1 then continue; fi;
    nn := nn + 1;
    vals := List(stabs, S -> lpart(Size(Intersection(S, N)), l));
    m := Minimum(vals);
    if m^2 > nl then
      nviol := nviol + 1;
      Out("VIOLATION ", tag, " l=", l, " |N|_l=", nl, " min=", m);
    elif m^2 = nl then
      ntight := ntight + 1;
      mins := Filtered([1..Length(reps)], i -> vals[i] = m);
      # number of N-orbits of minimizers
      cnt := Sum(mins, i -> Index(H, ClosureGroup(N, stabs[i])));
      if cnt = 1 then
        nrigid := nrigid + 1;
        Out("RIGID ", tag, " l=", l, " |N|_l=", nl, " min=", m);
      fi;
    fi;
  od;
  Out("FULL ", tag, " l=", l, " |G|=", Size(H), " orbits=", Length(orbs),
      " normals=", nn, " tight=", ntight, " viol=", nviol,
      " rigid=", nrigid);
end;;

# treats the group IrreducibleSolubleMatrixGroup(n, q, d, k)
# if its order is at most cap
CheckGroup := function(n, q, d, k, cap)
  local FF, p, F, G, s, iso, H, pcgs, mats, pr, l, ok, i, v, S, tag;
  p := SmallestRootInt(q); F := GF(q);
  G := IrreducibleSolubleMatrixGroup(n, q, d, k); s := Size(G);
  if s > cap then return; fi;
  pr := Filtered(PrimeDivisors(s), x -> x <> p);
  if pr = [] then return; fi;
  # necessary conditions for case 28 of [Holt-Yang, Table 4.1],
  # see the proof of Proposition 6.9
  if n = 8 and q = 5 then
    if s mod 96 <> 0 then return; fi;
    if d = 1 and s mod 192 <> 0 then return; fi;
    if 552960 mod s <> 0 then return; fi;
  fi;
  iso := IsomorphismPcGroup(G); H := Image(iso);
  if n = 8 and q = 5 then
    FF := FittingSubgroup(H);
    if IsAbelian(FF) or Size(FF) mod 32 <> 0
       or not IsCyclic(Centre(FF)) then return; fi;
  fi;
  for l in pr do
    # (G,V) must lie in X_l
    if not IsPNilpotent(H, l) or Size(PCore(H, l)) = 1 then continue; fi;
    # pairs covered by Lemma 6.5
    if IsPGroup(FittingSubgroup(H))
       and PrimePGroup(FittingSubgroup(H)) = l then continue; fi;
    pcgs := Pcgs(H);
    mats := List(pcgs, x -> PreImagesRepresentative(iso, x));
    ok := false;
    for i in [1..25] do
      v := Random(F^n);
      if IsZero(v) then continue; fi;
      S := Stabilizer(H, v, pcgs, mats, OnRight);
      if Size(S) mod l <> 0 then ok := true; break; fi;
    od;
    tag := Concatenation("(", String(n), ",", String(q), ",",
                         String(d), ",", String(k), ")");
    if ok then
      Out("RANDOK ", tag, " l=", l);
    else
      FullAnalysis(H, pcgs, mats, F, n, l, tag);
    fi;
  od;
end;;

# all groups in GL(n,q) with splitting field of degree d
# and order at most cap
CheckSegment := function(n, q, d, cap)
  local k;
  for k in IndicesIrreducibleSolubleMatrixGroups(n, q, d) do
    CheckGroup(n, q, d, k, cap);
  od;
  Out("SEGDONE n=", n, " q=", q, " d=", d, " cap=", cap);
end;;

# the capped runs for the cases of Table 1 and Proposition 6.16
for seg in [[4,11,1,11520], [4,13,1,13824],
            [6,5,1,10368], [6,5,2,10368], [6,5,3,2976],
            [6,7,1,8208], [6,7,3,8208],
            [12,2,1,12960], [12,2,2,12960], [12,2,4,12960],
            [8,3,1,7680], [8,3,2,7680], [8,3,4,7680],
            [8,5,1,55296], [8,5,2,55296]] do
  CheckSegment(seg[1], seg[2], seg[3], seg[4]);
od;

# the exhaustive runs of Remark 6.10
for nq in [[2,3],[2,5],[2,7],[2,11],[2,13],[2,17],[2,19],[2,23],
           [3,5],[3,7],[3,13],[3,19],[4,3],[4,5],[4,7],
           [5,2],[5,3],[6,2],[6,3]] do
  for d in DivisorsInt(nq[1]) do
    CheckSegment(nq[1], nq[2], d, infinity);
  od;
od;
\end{lstlisting}

The vectors found by the random search were not stored, but a new run of the code produces new such vectors. In our runs, the log files contained no line starting with \texttt{VIOLATION} or \texttt{RIGID}. Table \ref{table2} summarizes the results. Here ``groups'' is the number of irreducible solvable subgroups (up to conjugacy) in the given range, ``$p'$-vector'' is the number of pairs $(G,p)$ with $(G,V)\in\X_p$ for which a vector with $p'$-centralizer was found, and ``full'' is the number of pairs for which a complete orbit analysis was carried out. In the runs marked with $\dagger$, the pairs covered by Lemma \ref{lemmaF} were skipped and are not counted, and for $\GL(8,5)$ the count of groups includes the groups excluded by the filters described in the proof of Proposition \ref{nonmeta}.

\begin{table}[ht]
\centering
\small
\begin{tabular}{cccrrrr}
\toprule
$d$ & $r$ & splitting degrees & bound for $|G|$ & groups & $p'$-vector & full\\
\midrule
2 & 3 & all & -- & 7 & 3 & 2\\
2 & 5 & all & -- & 19 & 13 & 5\\
2 & 7 & all & -- & 29 & 24 & 7\\
2 & 11 & all & -- & 42 & 51 & 6\\
2 & 13 & all & -- & 62 & 64 & 10\\
2 & 17 & all & -- & 75 & 63 & 9\\
2 & 19 & all & -- & 77 & 93 & 9\\
2 & 23 & all & -- & 54 & 69 & 8\\
3 & 5 & all & -- & 22 & 7 & 0\\
3 & 7 & all & -- & 62 & 42 & 10\\
3 & 13 & all & -- & 136 & 96 & 9\\
3 & 19 & all & -- & 185 & 133 & 13\\
4 & 3 & all & -- & 108 & 30 & 39\\
4 & 5 & all & -- & 509 & 272 & 131\\
4 & 7 & all & -- & 894 & 682 & 167\\
5 & 2 & all & -- & 2 & 1 & 0\\
5 & 3 & all & -- & 16 & 8 & 0\\
6 & 2 & all & -- & 40 & 10 & 3\\
6 & 3 & all & -- & 324 & 69 & 16\\
\midrule
4 & 11 & 1 & 11520 & 938 & 824 & 121\\
4 & 13 & 1 & 13824 & 2692 & 2226 & 255\\
6 & 5 & 1, 2 & 10368 & 1290 & $249^\dagger$ & $12^\dagger$\\
6 & 5 & 3 & 2976 & 63 & $71^\dagger$ & $6^\dagger$\\
6 & 7 & 1, 3 & 8208 & 1979 & $885^\dagger$ & $19^\dagger$\\
8 & 3 & 1, 2, 4 & 7680 & 5643 & $641^\dagger$ & $413^\dagger$\\
8 & 5 & 1, 2 & 55296 & 98950 & $10593^\dagger$ & $16^\dagger$\\
12 & 2 & 1, 2, 4 & 12960 & 457 & $74^\dagger$ & $10^\dagger$\\
\bottomrule
\end{tabular}
\caption{Summary of the computations.}\label{table2}
\end{table}

In total, $1296$ pairs required a complete orbit analysis. In more than $240$ of them, some normal subgroup acted tightly, and in each of these cases the minimizers formed at least two orbits.

\begin{center}
    \textbf{Acknowledgement}
\end{center}
The author was supported by a grant from the Simons Foundation (SFI-MPS-TSM-00014153 to Thomas M. Keller). The author also wouldn like to thank Yong Yang for some discussions during the early stages of the project.

\begin{center}
    \textbf{AI statement}
\end{center}
This project was started about eight years ago. A working draft existed when AI systems became as powerful as they are today, and so the author used Claude (Opus 5.5) by Anthropic to finish the project. In particular, the author devised the project
goal and the proof strategy. Before using AI, he had Theorem A as the goal and knew the statement of Theorem C(a) and
had the reduction in Section 4, and Corollary D. He also had
a special case of Theorem B and some of the preliminary results in Section 2. The paper was then completed in an extended interaction between the author and Claude. The author posed the problem and provided the draft, and some of the preliminary results. The AI system found the general examples in Theorem B,, developed most of the arguments in Sections 5 and 6, carried out the computer calculations described in Section 6, and produced a first version of the manuscript, which was edited by the author. All arguments and computations were subsequently checked by the author, who takes full responsibility for the content of this paper.

\end{document}